\documentclass{article}

\usepackage{verbatim}
\usepackage{amsmath}
\usepackage{amsfonts}
\usepackage{amsthm}
\newtheorem{rem1}{Remark}[section]
\newtheorem{lem1}{Lemma}[section]
\newtheorem{cor1}{Corollary}[section]

\newtheorem{prop1}{Proposition}[section]
\newtheorem{thm1}{Theorem}[section]

\begin{document}
\title{Qualitative features of periodic solutions of KdV}\author{T. Kappeler\footnote{Supported in part by the Swiss National Science Foundation}, 
B. Schaad\footnote{Supported in part by the Swiss National Science Foundation}, P. Topalov\footnote{Supported in part by NSF DMS-0901443}}
\maketitle
\begin{abstract} \noindent 
In this paper we prove new qualitative features of solutions of KdV on the circle. The first result says that the Fourier
coefficients of a solution of KdV in Sobolev space $H^N,\, N\geq 0$, admit a WKB type expansion up to first order with
strongly oscillating phase factors defined in terms of the KdV frequencies. The second result provides estimates for the
approximation of such a solution by trigonometric polynomials of sufficiently large degree.

\noindent
Mathematics Subject Classification (2000): 35Q53, 35B40, 35B65, 37K10, 37K15, 37K40, 
\end{abstract}

\section{Introduction}\label{intro}
Consider the Korteweg-de Vries equation (KdV) 
\begin{equation}
\label{KdVequation}\partial_tu= -\partial^3_x u + 6u\partial_xu
\end{equation} 
on the circle $\mathbb{T}=\mathbb{R}/\mathbb{Z}$. It is globally in time well-posed on the Sobolev spaces
$H^s \equiv H^s(\mathbb{T},\mathbb{R})$ with $s\geq -1$ (\cite{KT}). 
The aim of this paper is to describe new qualitative features of periodic solutions of KdV.
First note that in contrast to solutions on the real line, periodic solutions do not have a special profile decomposition as 
$t\rightarrow \pm \infty $. Our main point of interest, related to the numerical experiments of Fermi, Pasta, and Ulam of 
particle chains, is to know how the distribution of energy among the Fourier modes evolves. 
A partial result in this direction says that due to the integrals provided by the KdV  hierarchy, the Sobolev norms of 
smooth solutions stay bounded uniformly in time. In this paper we make further contributions to the study of how the Fourier coefficients $\hat u_n(t)= \int_0^1 u(t,x)e^{-2\pi in x} dx$ of a solution $u(t,x)$ of \eqref{KdVequation} evolve in time. 
Our first result aims at describing dispersion phenomena for solutions of KdV by studying how  $\hat u_n(t)$ evolve for  $|n|$ large.
More precisely, we want to investigate if $\hat u_n(t)$ admits a WKB type expansion of the form
\begin{align}\label{intro1bis}\hat u_n(t)=e^{iw_nt}\left(a_n(t)+\frac{b_n(t)}{n} + \dots \right),\end{align}
where $e^{i w_nt}$ is a strongly oscillating phase factor with frequency $w_n$  and the coefficients $a_n(t),b_n(t),\dots$  vary more slowly and satisfy the estimates 
\[\sum n^{2N}|a_n(t)|^2< \infty \quad \text{and} \quad \sum n^{2N}|b_n(t)|^2< \infty.\]
To state our result more precisely, denote by $\omega_n,\,n\geq 1$, the KdV frequencies of $u(t)$. Let us recall how they are defined. The KdV equation can be written as a Hamiltonian PDE with phase space $L^2$ and Poisson bracket
\begin{align}\label{1bis}\{F,G\}(q):= \int_0^1 \partial F\,\partial_x \,\partial G dx \end{align}
where $F,G$ are $C^1$-functionals on $L^2$ and $\partial F$ denotes the $L^2$-gradient of $F$. Then KdV takes the form $\partial_t u=\partial_x\partial_u \mathcal{H}$ where $\mathcal{H}$ is the KdV Hamiltonian
\[
\mathcal{H}(q):=\int_0^1\left(\frac{1}{2}(\partial_xq)^2+q^3\right)dx.
\]
In terms of this set-up, the $\omega_n$'s are given by \[\omega_n= \partial_{I_n}\mathcal{H}.\]
Here we use that $\mathcal{H}$ can be expressed as a real analytic function of the action variables $I_n,\,n\geq 1$, so that the partial derivatives $\partial_{I_n}\mathcal{H}$ are well defined -- see below for more details. 
Alternatively, $\omega_n$ can be viewed as a function of $q$, which by a slight abuse of terminology, we also denote by $\omega_n$. Clearly, for any $n\geq 1,\, \omega_n(u(t))$ is independent of $t$ and depends in a nonlinear fashion on $u(0)$. 
It is convenient to introduce 
\[\omega_{-n}:=-\omega_n\quad\forall n\in \mathbb{Z}_{\geq 1} \quad  \text{and} \quad \omega_0:=0\]
and to denote the KdV flow by $S^t$, i.e., 
$S^t(u(0))= u(t)$. In addition, let  
\begin{equation}\label{eq:main_identity}
R^t(u(0)):= S^t(u(0))-\sum_{n\in \mathbb{Z}}  e^{i\omega_n t}\hat u_n(0)e^{2\pi i n x}
\end{equation}
where for any $n\in  \mathbb{Z}$, $\omega_n=\omega_n(u(0))$ .
\begin{thm1}\label{1smoothingthm2} For $q=u(0)\in H^N$, $N\in\mathbb{Z}_{\geq 0}$, the error 
$R^t(q)$ of the approximation  
$\sum_{n\in \mathbb{Z}}  e^{i\omega_n t} \hat q_n e^{2\pi i n x}$ of the flow $S^t(q)$ has the following properties:
\begin{itemize}\item[(i)] $R^t: H^N\to H^{N+1}$ is continuous;
\item[(ii)] for any $q\in H^N$, the orbit $\{R^t(q)|\, t\in \mathbb{R}\}$ is  relatively compact  in $H^{N+1}$;
\item[(iii)] for any $M>0$, the set of orbits $\{R^t(q) |\, t\in \mathbb{R},\, q\in H^N, \|q\|_{H^{N}}\leq M  \}$ is bounded in 
$H^{N+1}$;
\item[(iv)] if  in addition    $N\in \mathbb{Z}_{\geq 1}$, then  $\partial_t R^t: H^N \to H^{N-1}$ is continuous and for any $q\in H^N$, the orbit $\{\partial_t R^t(q)|\, t\in \mathbb{R}\}$ is relatively compact in $H^{N-1}$. 
Moreover for any $M>0$ the set of orbits $\{\partial_t R^t(q) |\, t\in \mathbb{R},\, q\in H^N, \|q\|_{H^{N}}\leq M  \}$ is bounded in $H^{N-1}$.
\end{itemize}
\end{thm1}

\begin{rem1}\label{top0.1}
Actually, one can prove that for any $c\in \mathbb{R}$, the restrictions of $R^t$ and $\partial_t R^t$ to the affine subspace $H_c^N=\{q\in H^N | \, \int_0^1q(x) dx =c\}$ are real analytic. 
See Remark \ref{toprem6.1} for a precise statement.
\end{rem1}

\begin{rem1} In case $u(0)$ is a finite gap potential, there are formulas, due to Its-Matveev \cite{IM}, for the frequencies $\omega_n$ in terms of periods of an Abelian differential, defined on the spectral curve associated to $u(0)$.
These formulas can be extended to potentials in $H^N,\, N\geq -1$, -- cf \cite{KT}. Alternative formulas can be found in \cite{KaP}, Appendix F. 
\end{rem1}

\begin{rem1} Note that the frequencies $\omega_n$ depend on the initial conditions in a nonlinear way.
The statements of Theorem \ref{1smoothingthm2} no longer hold if the KdV frequencies $\omega_n$ are replaced by their linearization at $0$, i.e., by $(2\pi n)^3$. 
Results on linear approximations of solutions of KdV were recently obtained by Erdogan and Tzirakis \cite{ETZ}, Theorem 1.2. 
In Appendix B, for initial data $q\in H^N$ with $N\in \mathbb{Z}_{\geq 1}$, we derive a stronger version of their result as a corollary of Theorem 1.1. 
\end{rem1}

\noindent In terms of the above WKB ansatz \eqref{intro1bis}, Theorem \ref{1smoothingthm2} says that with $w_n:= \omega_n$ and $a_n(t):= \hat u_n(0)$, the remainder term 
\[\rho_n(t):=b_n(t)+ \dots := n\cdot \left(e^{-i\omega_n t}\hat u_n(t)- \hat u_n(0)\right)= n\widehat R^t_n(u(0))e^{-i\omega_n t}\] satisfies $\sum n^{2N}|\rho_n(t)|^2< \infty$ and in case $N\in \mathbb{Z}_{\geq 1}$, 
\begin{align}\label{wkbosz}\sum n^{2(N-2)}|\partial_t\rho_n(t)|^2< \infty.\end{align} 
As the asymptotics of the KdV frequencies are given by 
$\omega_n= 8\pi^3 n^3 +  O(n)$ (cf formula \eqref{osz6ter} in Section \ref{section7bis}) estimate  \eqref{wkbosz} quantifies the assertion that $(\rho_n(t))_{n\in \mathbb{Z}}$ varies more slowly than $(\hat u_n(t))_{n\in \mathbb{Z}}$.

As an immediate consequence of Theorem \ref{1smoothingthm2} we obtain uniform bounds of the $H^s$-norms of 
the solutions of KdV in any fractional Sobolev space $H^s$, $s\ge 0$. Such bounds are of interest as the $H^s$-norms for 
$s\in\mathbb{R}_{\ge 0}\setminus\mathbb{Z}_{\ge 0}$ are not related to any Hamiltonian in the KdV-hierarchy.
To the best of our knowledge, they have {\em not} been known so far.
\begin{cor1}\label{coro:norms}
For any $s\in\mathbb{R}_{\ge 0}$ and for any $M>0$ the set
\[
\big\{S^t(q)\,\big|\,t\in\mathbb{R}, q\in H^s, \|q\|_{H^s}\le M\big\}
\]
is bounded in $H^s$.
\end{cor1}

\begin{rem1}
It is shown in \cite{CKSTT} that for any $-1/2\le s<0$ and for any $M>0$ and $\varepsilon>0$ there exists
$C>0$ such that for any $q\in H^s$ with $\|q\|_{H^s}\le M$ and for any $t\in\mathbb R$
\[
\|S^t(q)\|_{H^s}\le C (1+|t|)^{|s|+\varepsilon}\,.
\]
\end{rem1}

As a second application of Theorem \ref{1smoothingthm2} we obtain the following
\begin{cor1}\label{coro:kdv_weak}
For any $s\in\mathbb{R}_{\ge 0}$ and any $t\in\mathbb{R}$ the solution map $S^t : H^s\to H^s$ is weakly continuous, i.e.
for any sequence $(q_j)_{j\ge 1}$ in $H^s$ that converges weakly to $q$ in $H^s$ the sequence $(S^t(q_j))_{j\ge 1}$ converges weakly to $S^t(q)$ in $H^s$.
\end{cor1}
\noindent We note that recently such a result was proved for the NLS equation in $L^2$ by Oh and Sulem \cite{OS} by completely different methods.

\noindent The second result of this paper concerns the approximation of KdV solutions by trigonometric  polynomials. 
For any $L\in \mathbb{Z}_{\geq 1}$, denote by $P_L:L^2\to L^2$ the $L^2$-orthogonal projection of $L^2=H^0(\mathbb{T},\mathbb{R})$ onto the $2L+1$ dimensional 
$\mathbb{R}$-vector space generated by $e^{2\pi in x},\; |n|\leq L$.
\begin{thm1}\label{applthm5} Let $s\in\mathbb{R}_{\geq 0}$ be arbitrary. Then for any $M>0$ and $\varepsilon >0 $ there exists $L_{\varepsilon,M}\geq 1$ such that for any $u(0)\in H^s$, with 
$\|u(0)\|_{H^s}\leq M$, $L\geq L_{\varepsilon,M}$ and any $t\in\mathbb{R}$
\[ \|(Id-P_L)u(0)\|_{H^s} -\varepsilon \leq  \|(Id-P_L)u(t)\|_{H^s}\leq \|(Id-P_L)u(0)\|_{H^s} +\varepsilon.\]
In particular, if $u(0)$ with $\|u(0)\|_{H^s} \leq M$  is a trigonometric polynomial of order $L_*$, then for any 
$L\geq \max(L_*,L_{\varepsilon, M})$, $P_L u(t)$ approximates $u(t)$ \textit{uniformly in $t\in \mathbb{R}$} up to
an error of size $\varepsilon$.
\end{thm1}
\begin{rem1} 
The proof of Theorem \ref{applthm5} shows that for any $|n|> L_{\varepsilon, M}$ and $\|u(0)\|_{H^s}\leq M$,
\[|\hat u_n(0)|-\varepsilon \leq |\hat u_n(t)|\leq |\hat u_n(0)|+\varepsilon \quad \forall t\in \mathbb{R}.\]
It means that for $|n|$ sufficiently large, the amplitude of the $n$'th Fourier mode is approximately constant, uniformly on bounded sets of $H^s$.
\end{rem1}
\begin{rem1}
It follows from the proof of Theorem \ref{1smoothingthm2} that corresponding results hold for the flow of any Hamiltonian in the Poisson algebra of KdV. In particular, this is true for the flows of Hamiltonians in the KdV hierarchy.
\end{rem1}
The main ingredient of the proofs of Theorem \ref{1smoothingthm2} and Theorem \ref{applthm5} are refined asymptotics of the Birkhoff map of KdV.
This map provides normal coordinates,  allowing  to solve KdV by quadrature. Let us recall its set-up. First note that the average of any solution $u(t)\equiv u(t,x)$ of KdV in $H^s$ is a conserved quantity. 
In particular, for any $c\in\mathbb{R}$, KdV leaves the subspaces $H^s_c \equiv H^s_c(\mathbb{T},\mathbb{R})$ of $H^s$  invariant where 
\[
H^s_c = \left\{p(x)= \sum \hat p_n e^{2\pi inx}\;\big| \,\hat p_0=c,\, \left\|p\right\|_{H^s} < \infty, \, \hat p_{-n}= \overline{\hat p_{n}}\;\forall n\in \mathbb{Z}\right\}
\]
with
\[
\|p\|_{H^s}:=\Bigl(\sum |n|^{2s}|\hat p_{n}|^2 \Bigr)^{\frac{1}{2}}.
\] 
In the case $s=0$, we often write $L^2_c$ for $H^0_c$ and $\|p\|$ instead of $\|p\|_{H^0}$.
To describe the normal coordinates of KdV, let us introduce for any $\alpha \in \mathbb{R}$ the $\mathbb{R}$-subspace $\mathfrak{h}^\alpha$ of $\ell^{2,\alpha}$, given by
\[\mathfrak{h}^\alpha:=\left\{z= (z_n)_{n \neq 0}\in \ell^{2,\alpha}|\,z_{-n}=\overline z_n \forall n\geq 1\right\} \] where  
\[\ell^{2,\alpha}\equiv \ell^{2,\alpha}(\mathbb{Z}_{0},\mathbb{C}):=\left\{z= (z_n)_{n \neq 0}|\;\left\|z\right\|_\alpha < \infty \right\}, \]
 $\mathbb{Z}_0:= \mathbb{Z}\setminus \{ 0 \}$,
and 
\[\left\|z\right\|_\alpha:=\Bigl(\sum_{n\neq 0}|n|^{2\alpha}|z_n|^2\Bigr)^{\frac{1}{2}} .\]
The space $\mathfrak{h}^\alpha$ is endowed  with the standard Poisson bracket for which 
\[
\{z_n,z_{-n}\}=-\{z_{-n},z_n\}=2i
\] 
for any $n\geq 1$ whereas all other brackets between coordinate functions vanish. 
Furthermore we denote by $H^s_{0,\mathbb{C}}\equiv H^s_0(\mathbb{T},\mathbb{C}), \;L^2_{0,\mathbb{C}}\equiv L^2_0(\mathbb{T},\mathbb{C})$ and $\mathfrak{h}^\alpha_\mathbb{C}$  
the complexification of the spaces $H^s_0,\; L^2_0$, and $\mathfrak{h}^\alpha$. Note that $\mathfrak{h}^\alpha_\mathbb{C}= \ell^{2,\alpha}(\mathbb{Z}_0,\mathbb{C})$.
A detailed proof of the following result can be found in  \cite{KaP} -- cf also \cite{KaMa}.
\begin{thm1}\label{KdVtheorem}
There exist an open neighbourhood $W$ of $L^2_0$ in $L^2_{0,\mathbb{C}}$ and a real analytic map $\Phi: W \rightarrow \mathfrak{h}^{1/2}_\mathbb{C}$ with the following properties:
\begin{itemize}\item[(BC1)] For any $s\in \mathbb{R}_{\ge 0}$, the restriction of $\Phi$ to $H_0^s$ is a canonical, bianalytic diffeomorphism onto $\mathfrak{h}^{s+1/2}$.
\item[(BC2)] When expressed in the new coordinates, the KdV-Hamiltonian $\mathcal{H}\circ \Phi^{-1}$, defined on $\mathfrak{h}^{3/2}$, is a real analytic function of the action variables $I_n=(z_nz_{-n})/2$, $n\geq 1$ alone.
\item[(BC3)] The differential $\Phi_0\equiv d_0\Phi$ of $\Phi$ at $0$ is the weighted Fourier transform,
\begin{align}\label{BC32bis} \Phi_0(h)= \left(\frac{1}{\sqrt{|n|\pi}}\hat h_n\right)_{n\neq 0}
\end{align}
\end{itemize}
\end{thm1}
\noindent The coordinates  $z_n,\, n\neq 0$, are referred to as (complex) Birkhoff coordinates whereas $\Phi$ is called Birkhoff map. Note that in \cite{KaP} the Birkhoff map
is defined slightly differently by setting $\Phi(q)$ to be $(x_n,y_n)_{n\geq 1}$ where $x_n=(z_n+z_{-n})/2$ and $y_n=i(z_n-z_{-n})/2$. The fact that  KdV admits globally defined Birkhoff coordinates 
is a very special feature of KdV. In more physical terms it says that KdV, when considered with periodic boundary conditions, is a system of infinitely many coupled oscillators.
\begin{rem1}
A result similar to the one of Theorem \ref{KdVtheorem} holds for the defocusing NLS equation.  A detailed proof can be found in \cite{GKP}. Cf also \cite{MKV}.
\end{rem1}
The key ingredient of the proofs of Theorem \ref{1smoothingthm2} and Theorem \ref{applthm5} is the following result on the asymptotics of the Birkhoff map, which has an interest in its own.
\begin{thm1}\label{1smoothingthm} For $N\in\mathbb{Z}_{\geq 0}$, there exists an open neighbourhood $W_N$ of $H^N_0$ in $H^N_{0,\mathbb{C}}\cap W$ 
so that   $\Phi-\Phi_0$ maps $W_N$ into $\mathfrak{h}^{N+3/2}_\mathbb{C}$ and, as a map from $W_N$ to $ \mathfrak{h}^{N+3/2}_\mathbb{C}$, is analytic. 
Here $W$ is the neighbourhood of $L^2_0$ in $L^2_{0,\mathbb{C}}$ of Theorem 
\ref{KdVtheorem}.
Furthermore,  the restriction $A:=(\Phi-\Phi_0)_{| H^N_0}: H^N_0\to\mathfrak{h}^{N+3/2}$ is a bounded map, i.e. it is bounded on bounded subsets of $H^N_0$.
\end{thm1}
As an immediate application of Theorem \ref{1smoothingthm} we get the following 
\begin{cor1}\label{coro:bnf_weak}
For any $N\in\mathbb{Z}_{\ge 0}$, $\Phi : H^N_0\to\mathfrak{h}^{N+1/2}$ is a weakly continuous map, i.e. for any sequence $(q_j)_{j\ge 1}$ in $H^N_0$ that converges weakly to $q$ in $H^N_0$
the sequence  $(\Phi(q_j))_{j\ge 1}$ converges weakly to $\Phi(q)$ in $\mathfrak{h}^{N+1/2}$.
A corresponding result holds for $\Phi^{-1}$.
\end{cor1}

\begin{rem1} 
In \cite{Schaad} a result similar to the one stated in Theorem \ref{1smoothingthm} is proved  for the Birkhoff map of KdV constructed in \cite{KM}, where  the phase space is endowed with the Poisson bracket introduced by Magri. 
As an application, a corresponding result is then derived for the modified Korteweg-de Vries equation (mKdV) on $H^N$ with $N\geq 1$. Indeed, it was shown in \cite{KST2} that the Miura map $f\mapsto f'+f^2$ canonically embeds the symplectic leaves 
of the phase space of mKdV, endowed with the Poisson bracket \eqref{1bis}, into the phase space of KdV, endowed with the Magri bracket. (For a detailed study of the Miura map see \cite{KT1}.)
As a consequence, results similar to the ones of Theorem \ref{1smoothingthm2} and Theorem \ref{applthm5} can be proved for mKdV -- see \cite{Schaad}.
\end{rem1}
\begin{rem1} We expect that similar results as the ones of Theorem \ref{1smoothingthm}  can be proved for the defocusing NLS equation. As a consequence, results similar to the ones of Theorem \ref{1smoothingthm2} and 
Theorem \ref{applthm5} are expected to hold for this equation.
\end{rem1}
\begin{rem1}The asymptotic estimates obtained to prove Theorem \ref{1smoothingthm} can be used to derive a formula of the differential of the Birkhoff map $\Phi$ at $q=0$ by a short calculation (cf Appendix A). 
In addition they can be used to get a short proof of the Fredholm property  of the differential of $\Phi$ at any $q\in H^N_0$ (Corollary \ref{smoothcor2}).
\end{rem1}
\begin{rem1} Normalizing transformations such as the Birkhoff map are often viewed as nonlinear versions of the Fourier transform. In the case of KdV, Theorem \ref{1smoothingthm} provides a qualitative statement in this respect, 
saying that $\Phi$ is a \textit{weakly nonlinear} perturbation of the (weighted) Fourier transform.
\end{rem1}

\vspace{.5cm}

\noindent {\em  Related results:}
Recently, Kuksin and Piatnitski initiated a study of random perturbations with damping of the KdV equation \cite{KuPi}, \cite{Ku}. More precisely they are interested, 
how the KdV-action variables evolve under certain perturbed equations. For this purpose they express the perturbed KdV  equation in normal coordinates. Up to highest order, it is a linear differential equation if the nonlinear part 
$\Phi-\Phi_0$ of the Birkhoff map is $1$-smoothing, i.e. if it maps $H^N_0$ to $\mathfrak{h}^{N+3/2}$ for any $N\geq 0$. In their recent paper, Kuksin and Perelman \cite{KP} succeeded in showing that on a neighbourhood $U$ of 
the equilibrium point $q=0$, there exists a canonical, real analytic diffeomorphism $\Psi: U\rightarrow V$ with $V \subseteq \mathfrak{h}^{1/2}$ a neighbourhood of $0$ in $\mathfrak{h}^{1/2}$ providing Birkhoff coordinates for KdV so that 
$\Psi -\Psi_0$  is $1$-smoothing where $\Psi_0$ denotes the linearization of $\Psi$ at $q=0$ and coincides with $\Phi_0$.
They obtain the map $\Psi$ by generalizing  Eliasson's construction of a Birkhoff map near an equilibrium point of a finite dimensional integrable system to a class of integrable PDEs including the KdV equation. 
In order to apply Eliasson's construction, Kuksin and Perelman need coordinates for the KdV equation, provided in \cite{Ka}, as a starting point.  Eliasson's construction is based on Moser's path-method and, in general, 
cannot be extended to get global coordinates. However, for the study of random perturbations of KdV in \cite{Ku},  global Birkhoff coordinates for KdV are needed. In \cite{KP}, it was conjectured that there exists a globally defined Birkhoff map 
$\Psi$ so that $\Psi-\Psi_0$ is 1-smoothing. Note that Birkhoff maps are not uniquely determined. Theorem \ref{1smoothingthm} confirms that this conjecture holds true and that $\Psi$ can be chosen to be the Birkhoff map of Theorem \ref{KdVtheorem}  

\noindent The paper is organized as follows. In Section 2 we review asymptotic estimates of various spectral quantities, obtained in \cite{KST1}. In Section 3 and Section 4, 
these estimates are used to improve on asymptotic estimates of actions, angles, and Birkhoff coordinates, obtained in \cite{KaP}. In Section 5 we show Theorem \ref{1smoothingthm} and Corollary \ref{coro:bnf_weak}.
Finally, in Section 6, Theorem \ref{1smoothingthm2}, Corollary \ref{coro:kdv_weak}, and Theorem \ref{applthm5} are proved.

\noindent  For the convenience of the reader we now recall the ones most frequently used in this paper.
 For $q$ in $L^2_{0,\mathbb{C}}$, the Schr\"odinger operator $L_q:=-d_x^2+q$, considered on the interval $[0,2]$ with periodic boundary conditions, has a discrete spectrum, 
consisting of a sequence of complex numbers bounded from below. We list them lexicographically and with algebraic multiplicities,
\[\lambda_0\preceq \lambda_1\preceq \lambda_2 \preceq\lambda_3\preceq \dots\; \] where two  complex numbers $a,b,$ are ordered lexicographically, $a\preceq b$,  if $[\operatorname{Re}a< \operatorname{Re}b]$ or 
$[\operatorname{Re}a=\operatorname{Re}b\;\text{and}\,\operatorname{Im}a\leq\operatorname{Im}b]$. These eigenvalues satisfy the asymptotics
\[(\lambda_{2n}-n^2\pi^2)_{n\geq 1},(\lambda_{2n-1}-n^2\pi^2)_{n\geq 1} \in \ell^2\] or, expressed in a more convenient form,
 \[\lambda_{2n-1},\lambda_{2n}= n^2\pi^2 +\ell^2_n,\]
valid uniformly on bounded subsets of $L^2_{0,\mathbb{C}}$. In particular, this means that for any $R>0$ there exists $r>0$ so that for any $q \in L^2_{0,\mathbb{C}}$ with $\|q\|\leq R$, 
\begin{align}\label{1ter} |\lambda_k-n^2\pi^2|\leq r\pi^2 \quad \forall k\in \{2n,2n-1\},\,\forall n\geq 1.\end{align}  
For real $q$, the periodic eigenvalues are real and satisfy
\[\lambda_0< \lambda_1\leq \lambda_2 <\lambda_3\leq \lambda_4 <\dots\; .\]
Restricting ourselves to a sufficiently small neigbourhood $W$ of $L^2_0$ in $L^2_{0,\mathbb{C}}$, we can always ensure that the closed intervals 
\[G_n= \{(1-t)\lambda_{2n-1}+ t \lambda_{2n}|\;0\leq t\leq 1\},\quad n\geq 1,\] as well as
\[G_0= \{ t+ \lambda_{0}|\;-\infty < t\leq 0\}\] are disjoint from each other. By a slight abuse of terminology, for any $n\geq 1$, we refer to the closed interval $G_n$ as the $n$'th gap and to $\gamma_n:= \lambda_{2n}-\lambda_{2n-1}$, 
as the $n$'th gap length. We denote by $\tau_n$ the middle point of $G_n$, $\tau_n=(\lambda_{2n}+\lambda_{2n-1})/2$. Due to the asymptotic behaviour of the periodic eigenvalues, the $G_n$'s admit mutually disjoint neighbourhoods 
$U_n\subseteq \mathbb{C}$ with $G_n \subseteq U_n$ called \textit{isolating neighbourhoods}. Moreover, inside each $U_n$, we choose a circuit $\Gamma_n$ around $G_n$ with counterclockwise orientation. 
Both $U_n$ and $\Gamma_n$ can be chosen to be locally independent of $q$. For $q$ in a sufficiently small neighbourhood $W$ of $L^2_0$ in $L^2_{0,\mathbb{C}}$, the $U_n$'s with $n\geq n_0, \, n_0=n_0(q)$ sufficiently large, can be chosen
to be discs,  $U_n=\{\lambda \in \mathbb{C}|\;|\lambda-n^2\pi^2|<r\pi^2\}$, where $0< r\leq n_0$ and $n_0, r$ are chosen
so large that \eqref{1ter} holds. 
Such neighbourhoods will be called isolating neighbourhoods with parameters $n_0\geq 1$ and $ r>0$. In the course of this paper, $W$  will be shrunk several times, but we continue to denote it by $W$. 

\noindent By $\Delta(\lambda)\equiv  \Delta(\lambda,q)$ we denote the discriminant of $-d_x^2+q$, 
\[\Delta(\lambda)= \operatorname{tr} M(1,\lambda)\] where $M(x,\lambda)$ is the $2\times 2$ matrix whose columns $(y_i(x,\lambda),y_i'(x,\lambda))^{T},\, i=1,2$, are solutions of $-y''+q y=\lambda y$ 
with $M(0,\lambda)=Id_{2\times 2}$. The function $\Delta(\lambda)$ is entire and $\Delta^2(\lambda)-4$ has a product representation 
(see \cite{KaP}, Proposition B.10) \begin{align}\Delta^2(\lambda)-4 = 4 (\lambda_0- \lambda) 
\prod_{k\geq 1} \frac{(\lambda_{2k}-\lambda)
(\lambda_{2k-1}-\lambda)}{\pi_k^4}\end{align} 
where $\pi_k=k\pi$ for any $k\geq 1$. For $q$ in $L^2_{0,\mathbb{C}}$, 
we also need to consider the operator $-d_x^2 +q$ on $[0,1]$ with Dirichlet or Neumann boundary conditions. The corresponding spectra are again discrete, consisting of sequences of complex numbers, bounded from below. 
They are referred to as Dirichlet, respectively Neumann eigenvalues. We list them lexicographically and with their algebraic multiplicities \[ \mu_1\preceq \mu_2\preceq \mu_3 \dots \quad \text{and} \quad \eta_0\preceq \eta_1 \preceq \eta_2 \dots \;.\] 
They satisfy the asymptotics\[\mu_n,\eta_n=n^2\pi^2+\ell^2_n,\] valid uniformly on bounded subsets of $L^2_{0,\mathbb{C}}$. For real $q$, the Dirichlet and the Neumann eigenvalues are real and satisfy
\[\lambda_1\leq\mu_1\leq \lambda_2<\lambda_3\leq \mu_2\leq\lambda_4 < \dots \quad \text{and} \quad \eta_0 \leq \lambda_0 < \lambda_1 \leq \eta_1  \leq \lambda_2< \dots \;  .\]
 Restricting ourselves to a sufficiently small neighbourhood $W$ of $L^2_0$ in $L^2_{0,\mathbb{C}}$, we can assure that for any $q\in W$ there exist isolating neighbourhoods
 $U_n \subset \mathbb{C}$ so that $\mu_n,\eta_n\in U_n,\,n\geq 1$, whereas $\lambda_0$ and $\eta_0$ are not contained in any of the $U_n$'s. Isolating neighbourhoods with this additional property 
can be chosen to be locally independent of $q$. 

Finally let us recall the notion of the s-root, introduced in \cite{KaP}. 
For $a,b\in\mathbb{C}$, we define on $\mathbb{C}\setminus \{(1-t)a+t b\,|\;0\leq t\leq 1\}$ the s-root of $(b-\lambda)(\lambda-a)$, 
determined by setting for $\lambda \in \mathbb{C}$ with $|\lambda-\tau|>|b-a|$\[\sqrt[s]{(b-\lambda)(\lambda-a)}=i (\lambda-\tau)\sqrt[+]{1-w^2}\] where $\tau =(b+a)/2$ and 
$w=(b-a)/2(\lambda-\tau)$ -- see figure 2, p 62 in \cite{KaP}, showing a sign table. Here $\sqrt[+]\lambda$ denotes the principal branch of the square root on $\mathbb{C}\setminus (-\infty,0]$ 
characterized by, \[\sqrt[+]{\lambda}>0\quad \textrm{for}\; \lambda >0.\]
Throughout the paper, $\log \lambda$ denotes the principal branch of the logarithm, defined on $\mathbb{C}\setminus (-\infty,0]$. In particular, $\log 1=0$.

\vspace{0.5cm}

\noindent {\em Acknowledgment:} It is a pleasure to thank Wilhelm Schlag for useful discussions concerning Appendix B.

\section{Prerequisites}\label{section3}
In this section we review the asymptotic estimates of various spectral quantities, established in \cite{KST1} which are needed for the proof of Theorem \ref{1smoothingthm}. The results concern the asymptotics of the Floquet exponents 
$(\kappa_n)_{n\geq 1}$, the Dirichlet eigenvalues $(\mu_n)_{n\geq 1}$, the Neumann eigenvalues $(\eta_n)_{n\geq 0}$, and the periodic eigenvalues $(\lambda_n)_{n\geq 0}$ of the Schr\"odinger operator $-d_x^2 +q$ for a potential 
$q$ in $H^N_{0,\mathbb{C}}$ as well as the asymptotics of  $\gamma_n^2=(\lambda_{2n}-\lambda_{2n-1})^2$ and $\tau_n =(\lambda_{2n}+\lambda_{2n-1})/2,\;n\geq1$. 
Recall that for $q$ in $L^2_{0,\mathbb{C}}$,
\[\mu_n=n^2\pi^2+\ell^2_n,\quad \eta_n =n^2\pi^2 + \ell^2_n \quad\text{and}\quad \lambda_{2n},\lambda_{2n-1}= n^2\pi^2+ \ell^2_n.\]
For $f,g$ in $L^2_\mathbb{C}$, let \[\langle f,g\rangle= \int_0^1f(x)\bar g(x)dx.\]In particular, $\langle q, e^{2\pi i k x}\rangle = \int_0^1 q(x) e^{-2\pi i k x}dx$ denotes the k'th Fourier coefficient of $q$. The following asymptotics are known 
among experts in the field. 
\begin{thm1}\label{specthm} Let $q$ be in $H^N_{0,\mathbb{C}}$ with $N\in \mathbb{Z}_{\geq 0}$. Then 
\begin{align}\label{16a}\mu_n= m_n -\langle q,\cos2\pi nx\rangle +\frac{1}{n^{N+1}}\ell^2_n \\\label{n15bis}
\eta_n =m_n +\langle q,\cos2\pi nx\rangle+\frac{1}{n^{N+1}}\ell^2_n \end{align} uniformly on bounded subsets of potentials in $H^N_{0,\mathbb{C}}$. The quantity $m_n$ is of the form 
\begin{align}\label{12a} m_n= n^2\pi^2 +\sum_{2\leq 2j\leq N+1 }c_{2j}\frac{1}{n^{2j}}\end{align}
with coefficients $c_{2j}$ which are independent of $n$ and $N$ and given by integrals of polynomials in $q$ and its derivatives up to order $2j-2$. \end{thm1} 
\begin{rem1}The asymptotics \eqref{16a} are proved in \cite{Ma} and the uniform boundedness of the error in \eqref{16a} is shown in \cite{SS}.  A simple and self-contained proof of Theorem 2.1 can be found in \cite{KST1}.
\end{rem1}
\noindent Using similar arguments as in the proof of Theorem \ref{specthm}  corresponding estimates  for the periodic eigenvalues have been proved in \cite{KST1}.

\begin{thm1}\label{specthm2} Let $q$ be in $H^N_{0,\mathbb{C}}$ with $N\in \mathbb{Z}_{\geq 0}$. Then  \begin{align}\label{11a}\{\lambda_{2n},\lambda_{2n-1}\}= \{m_n \pm\sqrt{\hat q_n \hat q_{-n} +
\frac{1}{n^{2N+1}}\ell^2_n } +\frac{1}{n^{N+1}}\ell^2_n\} \end{align} uniformly on bounded subsets of potentials in $H^N_{0,\mathbb{C}}$. Again, $m_n$ is the expression defined in \eqref{12a}. 
\end{thm1}
\begin{rem1} A proof of the asymptotic estimate \eqref{11a}, but not of the uniform boundedness of the error terms, can be found in \cite{Ma}. 
\end{rem1}
\noindent Unfortunately, the asymptotics of Theorem \ref{specthm2} do not suffice for our purposes. Actually we need estimates of $\gamma_n^2= (\lambda_{2n}-\lambda_{2n-1})^2$ and $\tau_n= (\lambda_{2n}+
\lambda_{2n-1})/2$, $n\geq 1$ which are better than the ones obtained  from Theorem \ref{specthm2}.
\begin{thm1} \label{nthm3.2bis} Let $N\in \mathbb{Z}_{\geq 0}$. Then for any $q$ in $H^N_{0,\mathbb{C}}$, 
\begin{align}\label{n19} \gamma_n^2 = 4\hat q_n \hat q_{-n}  + \frac{1}{n^{2N+1}}\ell^1_n
\end{align}
uniformly on bounded subsets in $H^N_{0,\mathbb{C}}$.
\end{thm1}
\begin{proof} The claimed estimate follows from Theorem 1.2 in \cite{KaMi}. In the case at hand it says that 
\begin{align}\label{n17ter} \min_{\pm}\left|\gamma_n\pm 2\sqrt{\rho(n)\rho(-n)}\right|= \frac{1}{n^{N+1}}\ell_n^2
\end{align} uniformly on bounded sets of $H^N_{0,\mathbb{C}}$. Here the sequence $(\rho( n))_{n\in \mathbb{Z}}$ is given by 
\[\rho(n):= \langle q,e^{2\pi i n x}\rangle + \beta_1( n)\] with 
\[\beta_1( n):= \frac{1}{\pi^2} \sum_{k\neq \pm n} \frac{\langle q,e^{2\pi i (n-k) x}\rangle}{n-k}\frac{\langle q,e^{2\pi i (n+k) x}\rangle}{n+k}.\]
Then \begin{align}\nonumber & \left|\gamma_n^2 - 4\langle q, e^{2\pi i nx}\rangle\langle q, e^{-2\pi i nx}\rangle\right|\\\label{n17quater}\leq & \left|\gamma_n^2 -4\rho(n)\rho(-n)\right| + 
4\left| \rho(n)\rho(-n)-\langle q, e^{2\pi i nx}\rangle\langle q, e^{-2\pi i nx}\rangle \right|.
\end{align}
Note that \[\gamma_n^2 -4\rho(n)\rho(-n)= \left( \gamma_n + \varepsilon_n 2 \sqrt{\rho(n)\rho(-n)}\right)\left( \gamma_n -\varepsilon_n 2 \sqrt{\rho(n)\rho(-n)}\right),\]
where  $\varepsilon_n \in \{+,-\}$ is chosen such that 
\[\left|\gamma_n + \varepsilon_n2 \sqrt{\rho(n)\rho(-n)}\right| = \min_{\pm}\left|\gamma_n \pm 2\sqrt{\rho(n)\rho(-n)}\right|.\]
By Cauchy's inequality 
\begin{align*}
&\sum_{n\geq 1}n^{2N+1}|\gamma_n^2-4\rho(n)\rho(-n)| \\ \leq & \Big(\sum_{n\geq 1}n^{2N+2}|\gamma_n+\varepsilon_n 2\sqrt{\rho(n)\rho(-n)}|^2\Big)^{\!\!1/2}\!\!
\Big(\sum_{n\geq 1}n^{2N}\left|\gamma_n - \varepsilon_n 2 \sqrt{\rho(n)\rho(-n)}\right|^2\Big)^{\!\!1/2}\!\!\!.
\end{align*}
By \eqref{n17ter}, the first factor of the latter product is uniformly bounded on bounded sets of $q$'s in $W\cap H^N_{0,\mathbb{C}}$ where $W$ is an open neighborhood of the real space $H^N_0$ in $H^N_{0,\mathbb{C}}$, 
whereas the second factor can be estimated by 
\[
\left(\sum_{n\geq 1}n^{2N} |\gamma_n|^2\right)^{1/2}+ 2\left(\sum_{n\geq 1}n^{2N} |\rho(n)|^2\right)^{1/2}\left(\sum_{n\geq 1}n^{2N} |\rho(-n)|^2\right)^{1/2}.
\]
By \cite{KaMi}, Theorem 1.1, $\left(\sum_{n\geq 1}n^{2N} |\gamma_n|^2\right)^{1/2}$ is uniformly bounded on bounded sets of $q$'s in $W$ whereas by the definition of $\rho(\pm n)$ 
\[\left(\sum_{n\geq 1}n^{2N} |\rho(\pm n)|^2\right)^{1/2}\leq \|q\|_N + \|\beta_1\|_N \leq  \|q\|_N + \|q\|_N^2.\]
For the latter inequality we used that by \cite{KaMi}, Lemma 2.10, $\|\beta_1\|_{N+1}\leq \|q\|_N^2$.
It remains to estimate the second summand on the right hand side of \eqref{n17quater}. By the definition of $\rho(n)$
\begin{align*} & \sum_{n\geq 1}n^{2N+1} \left|\rho(n)\rho(-n) - \langle q, e^{2\pi i n x}\rangle\langle q, e^{-2\pi i n x}\rangle\right| \\\leq & 2\|\beta_1\|_{N+1} \|q\|_N + \|\beta_1\|_{N+1}^2 
\\ \leq & 2\|q\|_N^3(1+\|q\|_N),
\end{align*}
where we again used \cite{KaMi}, Lemma 2.10.
\end{proof}

For the sequences $(\tau_n)_{n\geq 1}$ the following asymptotic estimates are proved in \cite{KST1}.
\begin{thm1}\label{n3.2ter} (i) For any $q\in H^N_0,$ $N\in \mathbb{Z}_{\geq 0},$ 
\begin{align}\label{n17quinto} \tau_n(q)= m_n +\frac{1}{n^{N+1}}\ell^2_n
\end{align} where $m_n$ is given by \eqref{12a} and the error term is uniformly  bounded on bounded sets of potentials in $H^N_0$.

(ii) For any $N\in \mathbb{Z}_{\geq 0}$, there exists an open neighbourhood $W_N\subseteq  H^N_{0,\mathbb{C}}$ of $H^N_0$ so that \eqref{n17quinto} holds on $W_N$ with a locally uniformly bounded error term. 
\end{thm1}
\noindent Combining Theorem \ref{specthm} and  Theorem \ref{n3.2ter} one obtains
\begin{cor1}\label{corollary3.2}  (i) For any $q\in H^N_0,\, N\in \mathbb{Z}_{\geq 0}$, 
\begin{align}\label{17a}\tau_n-\mu_n&= \langle q,\cos2\pi nx\rangle +\frac{1}{n^{N+1}}\ell^2_n 
\end{align} 
where the error term is uniformly  bounded on bounded sets of potentials in $H^N_0$.

(ii) For any $N\in \mathbb{Z}_{\geq 0}$, there exists an open neighbourhood $W_N\subseteq  H^N_{0,\mathbb{C}}$ of $H^N_0$ so that \eqref{17a} holds on $W_N$ with a locally uniformly bounded error term.  
\end{cor1}

Furthermore we need asymptotic estimates for the Floquet exponents, defined by
\begin{align}\label{44bis.1}
\kappa_n= \log \left((-1)^ny_2'(1,\mu_n)\right).
\end{align}
 Here $\mu_n= \mu_n(q)$ is the n'th Dirichlet eigenvalue of $L(q)= -d_x^2+q$ with $q\in L^2_0$ and $y_2(x,\lambda)$ is the fundamental solution of $-y''+qy =\lambda y$ satisfying $y_2(0,\lambda)=0$ and $y_2'(0,\lambda)=1$. 
Note that $y_2'(1,\mu_n)\neq 0$. Actually, it turns out that $(-1)^ny_2'(1,\mu_n)>0$ for any $q\in L^2_0$. Hence $\log\left((-1)^ny_2'(1,\mu_n)\right)$ is well-defined with $\log$ denoting the principal branch of the logarithm. In fact, 
there exists a neighbourhood $W$ of $L^2_0$ in $L^2_{0,\mathbb{C}}$ so that for $q\in W$ and any $n\geq 1$, $\kappa_n(q)$ is well-defined by \eqref{44bis.1}. The $\kappa_n$'s have been introduced in \cite{FM} and studied 
for square integrable potentials in \cite{PT}. Note that for $\lambda= \mu_n$ the Floquet matrix 
\[M(\lambda)=\begin{pmatrix} y_1(1,\lambda) & y_2(1,\lambda) \\ y_1'(1,\lambda) & y_2'(1,\lambda)
\end{pmatrix}\]
is lower triangular. Hence $y_2'(1,\mu_n)$ is one of the two Floquet multipliers of $M(\mu_n)$, the other one being $y_1(1,\mu_n)$ which by the Wronskian identity equals $1/y_2'(1,\mu_n)$. In particular it follows that 
\begin{align}\label{B.1} \kappa_n=- \log \left((-1)^ny_1(1,\mu_n)\right). 
\end{align}
In \cite{KST1}  we prove 
\begin{thm1}\label{4bis.1} Let $N\geq 0$. Then for any $q$ in $W\cap H^N_{0,\mathbb{C}}$, \[\kappa_n= \frac{1}{2\pi n} \left(\langle q, \sin 2\pi n x\rangle+ \frac{1}{n^{N+1}}\ell^2_n\right)\]
uniformly on bounded subsets of $W\cap H^N_{0,\mathbb{C}}$.
\end{thm1}
Note that in contrast to the asymptotics of the Dirichlet eigenvalues or the Neumann eigenvalues, the size of $\kappa_n$ for any $q\in H^N_0$ is of the order of $\frac{1}{n^{N+1}}\ell^2_n$. The case $N=0$ is much simpler and has 
been treated in \cite{PT}, p 60.

\noindent Finally we state some applications of the asymptotics of the periodic eigenvalues. For our purposes it suffies to consider potentials $q$ in a sufficiently small neighbourhood $W$ of $L^2_0$ in $L^2_{0,\mathbb{C}}$.
\noindent Recall that we denote by $\Delta(\lambda)$ the discriminant of $-d^2_x+ q$ and that $\Delta^2(\lambda)-4$ has the product 
representation 
\[\Delta^2(\lambda)-4= 4(\lambda_0-\lambda)\prod_{n\geq 1}\frac{(\lambda_{2n}-\lambda)(\lambda_{2n-1}-\lambda)}{\pi_n^4}\] where $\pi_n=n\pi$ for any $n\geq 1$. Let $\dot \Delta(\lambda) = \partial_{\lambda}\Delta(\lambda)$.
According to Proposition B.13 of \cite{KaP} it also admits a product representation, \[\dot\Delta(\lambda) =-\prod_{n\geq 1}\frac{\dot\lambda_n-\lambda}{\pi_n^2},\]
and the zeros $\dot\lambda_n$ satisfy \begin{equation}\label{17bis}\dot\lambda_n -\tau_n= O(\gamma_n^2)\end{equation} locally uniformly for $q$ in $W$. Shrinking $W$, if necessary, we can assume without loss of generality that 
for any $n\geq 1,\; \dot\lambda_n\in U_n$ locally uniformly in $W$.  The following estimate improves on the one of Propostition B.13 in \cite{KaP}.
\begin{prop1}\label{proposition3.4}
For $q$ in $W$, \begin{equation}\label{3.7}\dot\lambda_n -\tau_n= \frac{\gamma_n^2}{n}\ell^2_n\end{equation} locally uniformly on $W$. On $L^2_0$, \eqref{3.7} is uniformly bounded on bounded subsets of $L^2_0$.
\end{prop1}

\begin{proof}For any given  $n\geq 1$, write \begin{align}\label{18bis}\Delta^2(\lambda)-4= \frac{(\lambda_{2n}-\lambda)(\lambda- \lambda_{2n-1})}{\pi_n^2}\Delta_n(\lambda)\end{align} where 
\begin{align}\label{19bis}\Delta_n(\lambda)= 4 \,\frac{\lambda -\lambda_0}{\pi_n^2}\left(\prod_{m\neq n}\frac{\lambda_{2m}-\lambda}{\pi_m^2}\right)\left(\prod_{m\neq n}\frac{\lambda_{2m-1}-\lambda}{\pi_m^2}\right).\end{align} 
Uniformly for  $\lambda\in U_n$,
\[\frac{\lambda-\lambda_0}{\pi_n^2}=1 + O\left(\frac{1}{n^2}\right)= 1+ \frac{1}{n} \ell^2_n.\]
By Corollary 7.1 in \cite{KST1}, uniformly for $\lambda\in U_n$,  
\begin{align*}\left(\prod_{m\neq n}\frac{\lambda_{2m}-\lambda}{\pi_m^2}\right)\left(\prod_{m\neq n}\frac{\lambda_{2m-1}-\lambda}{\pi_m^2}\right)&= \left(\frac{(-1)^{n+1}}{2}+ \frac{1}{n}\ell^2_n\right)^2\\&= 
\frac{1}{4}+ \frac{1}{n}\ell^2_n,
\end{align*}
hence \begin{align}\label{19ter} \Delta_n(\lambda)= 1 + \frac{1}{n}\ell^2_n \end{align} and by shrinking  the size of the isolating neighbourhoods one obtains from Cauchy's estimate that
\[\dot\Delta_n(\lambda)=  \frac{1}{n}\ell^2_n\] uniformly in $n\geq 1, \, \lambda \in U_n$. By \eqref{18bis}
\begin{align*}0&= \frac{d}{d\lambda}|_{\lambda=\dot\lambda_n}(\Delta^2(\lambda)-4)\\&= \frac{\lambda_{2n-1}-\dot\lambda_n + \lambda_{2n}- \dot\lambda_n}{\pi_n^2}\Delta_n(\dot \lambda_n)\\&+ 
\frac{(\lambda_{2n}-\dot\lambda_n)  (\dot\lambda_{n}- \lambda_{2n-1})}{\pi_n^2}\dot\Delta_n(\dot \lambda_n)
\end{align*}
or \[0= 2(\tau_n-\dot\lambda_n)(1+ \frac{1}{n}\ell^2_n)+ \left(\frac{\gamma_n^2}{4}- (\tau_n-\dot\lambda_n)^2\right)\frac{1}{n}\ell^2_n.\] 
Hence
\[(\tau_n-\dot\lambda_n)\left(1+ \frac{1}{n}\ell^2_n+(\tau_n-\dot\lambda_n)\frac{1}{n}\ell^2_n \right)= \frac{\gamma_n^2}{n}\ell^2_n.\]
As by \eqref{17bis}, $\tau_n-\dot\lambda_n= O(1)$  it then follows that \[\tau_n-\dot\lambda_n=\frac{\gamma_n^2}{n}\ell^2_n\] as claimed. Going through the arguments of the proof one sees that \eqref{3.7} holds locally uniformly 
on $W$ and uniformly on bounded subsets of $L^2_0$.
\end{proof}
\noindent By Theorem D.1 in \cite{KaP} there exists a sequence $(\psi_n)_{n\geq 1}$ of entire functions,
\begin{align}\label{21bis}\psi_n(\lambda)=\frac{2}{\pi_n}\prod_{m\neq n}\frac{\sigma^n_m-\lambda}{\pi_m^2} \end{align} so that
\begin{equation}\label{3.10} \frac{1}{2\pi}\int_{\Gamma_m}\frac{\psi_n(\lambda)}{\sqrt[c]{\Delta^2(\lambda)-4}}d\lambda=\delta_{mn}\quad \forall m\geq 1,
\end{equation}
where $\sqrt[c]{\Delta^2(\lambda)-4}$ denotes the canonical root introduced in \cite{KaP}, section 6.
Recall that $\Gamma_m$ denotes a counterclockwise oriented circuit in $U_m$ around the interval $G_m$ and $\tau_m=(\lambda_{2m}+\lambda_{2m-1})/2$.
By Theorem D.1 in \cite{KaP}, for any $n\geq 1,$ the zeros $\sigma_m^n,\, m\neq n,$ are real analytic functions of $q\in W$ so that $\sigma_m^n-\tau_m= O\left(\gamma_m^2/m\right)$ uniformly for
 $n\geq 1$, and locally uniformly in $q\in W$. Moreover one can choose $W$ so that $\sigma_m^n\in U_m$ 
for any $n\geq 1,\, m\neq n$, and locally uniformly in $W$.
\begin{prop1}\label{proposition3.5}For $q$ in $W$, \begin{align}\label{a24bis}\sigma^n_m-\tau_m= \frac{\gamma_m^2}{m}\ell^2_m \end{align}   locally uniformly on $W$ and uniformly in $n$.
On $L^2_0$, \eqref{a24bis}  is uniformly bounded on bounded subsets of $L^2_0$.
\end{prop1}

\begin{proof}  We drop the superscript in $\sigma^n=(\sigma_m^n)_{m\neq n}$ for the course of this proof. For $m \neq n$ one has
\begin{align}\label{top26bis}0=\int_{\Gamma_m}\frac{\sigma_m-\lambda}{\sqrt[s]{(\lambda_{2m}-\lambda)(\lambda- \lambda_{2m-1})}}\chi_m^n(\sigma,\lambda)d\lambda\end{align}
 with 
\[\chi_m^n(\sigma, \lambda)= (-1)^{m-1}\frac{\sigma_n -m^2\pi^2}{\sigma_n-\lambda}\frac{m\pi}{\sqrt[+]{\lambda -\lambda_0}}\prod_{j\neq m}\frac{\sigma_j-\lambda}{\sqrt[+]{(\lambda_{2j}-
\lambda)(\lambda_{2j-1}-\lambda)}}\] and $\sigma_n =\tau_n$.
For $\lambda=m^2\pi^2+\ell^2_m$, one has uniformly in $n\geq 1$, $m\neq n$, 
\begin{align*}\frac{\sigma_n-m^2\pi^2}{\sigma_n-\lambda}=1+\frac{\lambda-m^2\pi^2}{\sigma_n-\lambda}= 1+ \frac{1}{m}\ell^2_m
\end{align*} and $\sqrt[+]{\lambda-\lambda_0}=m\pi\sqrt[+]{1+\frac{1}{m}\ell^2_m}= m\pi\left(1+\frac{1}{m}\ell^2_m\right).$ 
Furthermore, by Proposition 7.1 in \cite{KST1},
\begin{equation}\label{19biss}\left(\prod_{j\neq m}\frac{\sigma_j-\lambda}{\sqrt[+]{(\lambda_{2j}-\lambda)(\lambda_{2j-1}-\lambda)}}\right)^2=1+ \frac{1}{m}\ell^2_m.\end{equation}
Taking square roots on both sides of \eqref{19biss}, one has \[(-1)^{m-1}\prod_{j\neq m}\frac{\sigma_j-\lambda}{\sqrt[+]{(\lambda_{2j}-\lambda)(\lambda_{2j-1}-\lambda)}}=1+ \frac{1}{m}\ell^2_m. \]
Altogether we have for $\lambda=m^2\pi^2 +\ell^2_m$, $\lambda \in U_m,$  \begin{equation}\label{omegadef}\chi_m^n(\lambda)=1+ \frac{1}{m}\ell^2_m\end{equation} uniformly in $n$.  
The integral 
\begin{align}\label{24bbis} 
\int_{\Gamma_m}\frac{\sigma_m-\lambda}{\sqrt[s]{(\lambda_{2m}-\lambda)(\lambda- \lambda_{2m-1})}}d\lambda 
\end{align} can be explicitly computed. 
In the case where $\lambda_{2m}=\lambda_{2m-1}$, one gets from the definition of the s-root and Cauchy's formula that the integral equals $2\pi(\sigma_m-\tau_m)$. In  the case $\lambda_{2m}\neq\lambda_{2m-1}$,
\[\int_{\Gamma_m}\frac{\sigma_m-\lambda}{\sqrt[s]{(\lambda_{2m}-\lambda)(\lambda- \lambda_{2m-1})}}d\lambda= 
2\int_{\lambda_{2m-1}}^{\lambda_{2m}}\frac{\sigma_m-\lambda}{\sqrt[+]{(\lambda_{2m}-\lambda)(\lambda- \lambda_{2m-1})}}d\lambda.\] By the change of coordinate $\lambda= \tau_m+ t\frac{\gamma_m}{2}$ we get
\begin{align*}\int_{\Gamma_m}\frac{\sigma_m-\lambda}{\sqrt[s]{(\lambda_{2m}-\lambda)(\lambda- \lambda_{2m-1})}}d\lambda&= 
2\int_{-1}^1\frac{\sigma_m-\tau_m}{\sqrt[+]{1-t^2}}d t -\gamma_m\int_{-1}^1\frac{t}{\sqrt[+]{1-t^2}}d t\\&=2 (\sigma_m- \tau_m) \int_{-1}^1\frac{1}{\sqrt[+]{1-t^2}}d t.
\end{align*} 
Therefore in both cases\[\frac{1}{2\pi }\int_{\Gamma_m}\frac{\sigma_m-\lambda}{\sqrt[s]{(\lambda_{2m}-\lambda)(\lambda- \lambda_{2m-1})}}d\lambda=\sigma_m-\tau_m.\] Hence, by \eqref{top26bis}  
\begin{align}\label{top29bis}0= (\sigma_m- \tau_m) \chi_m^n(\tau_m)+ \frac{1}{2\pi }\int_{\Gamma_m}
\frac{(\sigma_m-\lambda)(\chi_m^n(\lambda)-\chi_m^n(\tau_m))}{\sqrt[s]{(\lambda_{2m}-\lambda)(\lambda- \lambda_{2m-1})}}d\lambda.
\end{align} 
As $\tau_m= m^2\pi^2 +\ell^2_m$ it follows from \eqref{omegadef} that
\begin{equation}\label{20biss}\chi_m^n(\tau_m)= 1+ \frac{1}{m}\ell^2_m\end{equation}  and, by the Taylor expansion of $\chi^n_m$ at $\tau_m$ and Cauchy's estimate, for $\lambda =m^2\pi^2 +\ell^2_m$, $\lambda \in U_m$  
\[\frac{\chi_m^n(\lambda)-\chi_m^n(\tau_m)}{\lambda -\tau_m}= \frac{1}{m}\ell^2_m.\]
Finally as $\sigma_m -\lambda=O(\gamma_m), \lambda \in G_m$, one has $(\sigma_m-\lambda)(\lambda-\tau_m)=O(\gamma_m^2)$ for $\lambda \in G_m$.
Hence by Lemma M.1 in \cite{KaP}
\begin{align*}&\int_{\Gamma_m}\frac{(\sigma_m-\lambda)(\chi_m^n(\lambda)-\chi_m^n(\tau_m))}{\sqrt[s]{(\lambda_{2m}-\lambda)(\lambda- \lambda_{2m-1})}}d\lambda
\\=&2\int_{G_m}\frac{(\sigma_m-\lambda)(\lambda-\tau_m)\frac{\chi_m^n(\lambda)-\chi_m^n(\tau_m)}{\lambda-\tau_m}}{\sqrt[+]{(\lambda_{2m}-\lambda)(\lambda- \lambda_{2m-1})}}d\lambda \\=& \frac{\gamma_m^2}{m}\ell^2_m
\end{align*}
This together with \eqref{top29bis} gives \[ \sigma_m^n- \tau_m =
\frac{\gamma_m^2}{m}\ell^2_m \] uniformly in $n$. Going through the arguments of the proof one sees that \eqref{a24bis} holds locally uniformly on $W$ and uniformly on bounded subsets of $L^2_0$.
\end{proof}

\section{Asymptotics of actions and angles}

In this section we improve on estimates of the actions and angles obtained in 
\cite{KaP}, section 7 respectively section 8. Let us begin with asymptotics of the actions. Recall that $W$ is a (sufficiently small) neighbourhood of $L^2_0$ in $L^2_{0,\mathbb{C}}$. For $q$ in $W$, the n'th action variable is defined by 
\[I_n= \frac{1}{\pi}\int_{\Gamma_n} \lambda \frac{\dot\Delta(\lambda)}{\sqrt[c]{\Delta^2(\lambda)-4}}d\lambda \qquad (n\geq 1),
\]
where $\Delta(\lambda)$ is the discriminant, $\dot\Delta(\lambda)= \partial_\lambda\Delta(\lambda)$ and $\Gamma_n$ is a contour around the interval $G_n:=\{t\lambda_{2n-1}+(1-t)\lambda_{2n}|\;0\leq t\leq 1\}$ in
 the isolating neighbourhood $U_n$ of $G_n$ -- see at the end of the introduction or section 7 of \cite{KaP} for more details.
 First we need to derive improved estimates of the quotient
$I_n /\gamma_n^2$, given in \cite{KaP}, Theorem 7.3.
\begin{prop1}\label{quotientprop} Locally uniformly on $W$, the quotient $I_n /\gamma_n^2$ satisfies 
\begin{equation}\label{quotientestimate}
8\pi n \frac{I_n}{\gamma_n^2}= 1+\frac{1}{n}\ell^2_n.
\end{equation}  
Moreover
\begin{align}\label{15bis}
\xi_n = \sqrt[+]{8I_n /\gamma_n^2}= \frac{1}{\sqrt{\pi n}}(1+\frac{1}{n}\ell^2_n)
\end{align}  
is well-defined as a real analytic, non-vanishing function on W. In particular, at $q=0$, we have $\xi_n= \frac{1}{\sqrt{\pi n}}$ for all $n\geq 1$.
On $L^2_0$, \eqref{15bis} holds uniformly on bounded subsets of $L^2_0$.
\end{prop1}
\begin{proof}
We refer to \cite{KaP}, section 7,  for all notions, notations (and results) not explained here. In view of Theorem 7.3 in \cite{KaP} it only remains to be shown the improved asymptotics \eqref{quotientestimate}.  
Recall the product expansions
\begin{eqnarray}\nonumber \Delta^2(\lambda)-4&=& 4 (\lambda_0- \lambda) 
\prod_{n\geq 1} \frac{(\lambda_{2n}-\lambda)(\lambda_{2n-1}-\lambda)}{\pi_n^4}
\end{eqnarray} 
and
\begin{eqnarray}\nonumber
\dot\Delta(\lambda)&=& - \prod_{n\geq 1}\frac{\dot\lambda_n-\lambda}{\pi_n^2}.
\end{eqnarray}
For $\lambda$ on $\Gamma_n$ write
\begin{align}\label{27bbis}
\frac{\dot\Delta(\lambda)}{\sqrt[c]{\Delta^2(\lambda)-4}}=
\frac{1}{2\pi n}\frac{\lambda-\dot\lambda_n}{\sqrt[s]{(\lambda_{2n}-\lambda)(\lambda-\lambda_{2n-1})}}\chi _n(\lambda)
\end{align}
where for $\lambda\in U_n$,
\begin{align}
\label{4.0} \chi_n(\lambda)= (-1)^{n-1}\frac{n \pi}{\sqrt[+]{\lambda-\lambda_0}}
\prod_{m\neq n}\frac{\dot\lambda_m-\lambda}{\sqrt[+]{(\lambda_{2m}-\lambda)(\lambda_{2m-1}-\lambda)}}, 
\end{align}
and where the canonical root $\sqrt[c]{\Delta^2(\lambda)-4}$ has been introduced earlier. (For questions of signs, see section 6 in \cite{KaP}.)  
Note that for $\lambda\in U_n,$
\begin{align}\label{4.1}
\frac{n\pi}{\sqrt[+]{\lambda-\lambda_0}}= 1+\frac{1}{n}\ell^2_n.
\end{align}
Furthermore, by Proposition \ref{proposition3.4}, the roots $\dot\lambda_n$ of $\dot\Delta$ satisfy $\dot\lambda_n=\tau_n+ \frac{\gamma_n^2}{n}\ell^2_n$ and hence, by \cite{KST1}, Proposition 2.1, one has for $\lambda\in U_n$,
\begin{align*}
\left(\prod_{m\neq n}\frac{\dot\lambda_m-\lambda}{\sqrt[+]{(\lambda_{2m}-\lambda)(\lambda_{2m-1}-\lambda)}}\right)^2&= 
\prod_{m\neq n}\frac{\dot\lambda_m-\lambda}{\lambda_{2m}-\lambda}\prod_{m\neq n}\frac{\dot\lambda_m-\lambda}{\lambda_{2m-1}-\lambda}\\&= 1+\frac{1}{n}\ell^2_n
\end{align*}
or
\begin{align}\label{4.2}
(-1)^{n-1}\prod_{m\neq n}\frac{\dot\lambda_m-\lambda}{\sqrt[+]{(\lambda_{2m}-\lambda)(\lambda_{2m-1}-\lambda)}}=1+\frac{1}{n}\ell^2_n.
\end{align}
Combining \eqref{4.1} and \eqref{4.2} leads to the estimate 
\begin{equation}\label{4.3} 
\chi_n(\lambda) =1+ \frac{1}{n}\ell^2_n
\end{equation}
uniformly for $\lambda\in U_n$. Introduce \[ Z_n=\{q\in W|\; \lambda_{2n}(q)=\lambda_{2n-1}(q)\}.\]
Arguing as in section 7 of \cite{KaP} one gets for $q$ in $W\setminus Z_n$
\begin{align*}
 I_n & = \frac{1}{\pi}\int_{\Gamma_n} (\lambda-\dot\lambda_n) \frac{\dot\Delta(\lambda)}{\sqrt[c]{\Delta^2(\lambda)-4}}d\lambda\\&=
\frac{1}{2\pi^2n}\int_{\Gamma_n}\frac{(\lambda-\dot\lambda_n)^2}{\sqrt[s]{(\lambda_{2n}-\lambda)(\lambda-\lambda_{2n-1})}}\chi _n(\lambda)d\lambda.
\end{align*}
Substituting $\lambda= \tau_n+\zeta \gamma_n /2$ and setting 
$\delta_n= 2(\dot\lambda_n-\tau_n) /\gamma_n$ yields
\[I_n= \frac{\gamma_n^2}{8\pi^2 n}\int_{\Gamma_n'} (\zeta-\delta_n)^2 \chi_n(\tau_n+ \zeta\frac{\gamma_n}{2})\frac{d\zeta}{\sqrt[s]{1-\zeta^2}}\]
where $\Gamma_n'$ is some circuit in $\mathbb{C}$ around $[-1,1]$.
Thus on $W\setminus Z_n$, 
\begin{align*}8 \pi n\frac{I_n}{\gamma_n^2}&= \frac{1}{ \pi}\int_{\Gamma_n'} (\zeta-\delta_n)^2 \chi_n\left(\tau_n+ \zeta\frac{\gamma_n}{2}\right)\frac{d\zeta}{\sqrt[s]{1-\zeta^2}}\\
&=\frac{1}{ \pi}\int_{\Gamma_n'} (\zeta-\delta_n)^2 \chi_n(\tau_n)\frac{d\zeta}{\sqrt[s]{1-\zeta^2}}\\&+
\frac{1}{ \pi}\int_{\Gamma_n'} (\zeta-\delta_n)^2 \left(\chi_n\left(\tau_n+ \zeta\frac{\gamma_n}{2}\right)-\chi_n(\tau_n)\right)\frac{d\zeta}{\sqrt[s]{1-\zeta^2}}.
\end{align*}
By Proposition \ref{proposition3.4}, $\delta_n$ is well defined on all of $W$, hence so is the r.h.s. of the latter identity. As
\[(\zeta-\delta_n)^2= \zeta^2+\frac{\gamma_n}{n}\ell^2_n\]
Lemma M.1 in \cite{KaP} and \eqref{4.3}   then imply that 
\begin{align*}\frac{1}{\pi}\int_{\Gamma_n'}(\zeta-\delta_n)^2 \chi_n(\tau_n)\frac{d\zeta}{\sqrt[s]{1-\zeta^2}} &= \frac{1}{\pi}\int_{\Gamma_n'}\frac{\zeta^2 d\zeta}{\sqrt[s]{1-\zeta^2}} +\frac{1}{n}\ell^2_n\\&= 1+ \frac{1}{n}\ell^2_n.
\end{align*} By the Taylor expansion of order $0$ of $\chi_n$ at $\lambda=\tau_n$,  Cauchy's estimate, and \eqref{4.3} to bound $\dot\chi_n(\lambda)$, one gets
\begin{align*}\frac{1}{\pi}\int_{\Gamma_n'}(\zeta-\delta_n)^2 \left(\chi_n(\tau_n+\zeta\frac{\gamma_n}{2})-\chi_n(\tau_n)\right)\frac{d\zeta}{\sqrt[s]{1-\zeta^2}} = \frac{1}{n}\ell^2_n.
\end{align*}
Altogether,
\[8\pi n I_n/ \gamma_n^2= 1 +\frac{1}{n}\ell^2_n
\] and thus \[\xi_n = \frac{1}{\sqrt{\pi n}}(1+\frac{1}{n}\ell^2_n).\]
Going through the arguments of the proof one sees that the estimate \eqref{15bis} holds locally uniformly on $W$ and uniformly on bounded subsets of $L^2_0$.
\end{proof}
\vspace{0.4cm}

\noindent Proposition \ref{quotientprop} leads to the following asymptotics of the action variables. 
\begin{prop1}\label{actionasymptotics} Locally uniformly on $W\cap H^N_{0,\mathbb{C}}$
\begin{align}\label{a35bis} 2 I_n=\frac{1}{\pi n}\left(\hat q_n \hat q_{-n} +\frac{1}{n^{2N+1}}\ell^1_n\right).\end{align}
On $H^N_0$, the error  in \eqref{a35bis} is uniformly bounded on bounded subsets of $H^N_0$.
\end{prop1}
\begin{proof}
By Theorem \ref{nthm3.2bis}  one has
\[\gamma_n^2= 4 \hat q_n \hat q_{-n} +  \frac{1}{n^{2N+1}}\ell^1_n.
\]
Combined with  the asymptotics \eqref{quotientestimate} we then get,
\[2 I_n= \frac{1}{4n\pi}\gamma_n^2\left(1+\frac{1}{n}\ell^2_n\right)= \frac{1}{\pi n} (1+\frac{1}{n}\ell^2_n)\left(\hat q_n \hat q_{-n}
+\frac{1}{n^{2N+1}}\ell^1_n\right)
\]
or
\[2 I_n= \frac{1}{\pi n}\left(\hat q_n \hat q_{-n}  +\frac{1}{n^{2N+1}}\ell^1_n\right).\]
Going through the arguments of the proof one sees that \eqref{a35bis} holds locally uniformly on $W\cap H^N_{0,\mathbb{C}}$ and uniformly on bounded subsets of $H^N_0$.
\end{proof}
\noindent Next we improve on the estimates of the angle variables obtained in \cite{KaP}, section 8. To this end we use the improved estimate of the zeros $(\sigma_m^n)_{m\neq n}$ of the entire function  $\psi_n$,
\[\sigma_m^n= \tau_m+\frac{\gamma_m^2}{m}\ell^2_m\]
of Proposition \ref{proposition3.5}.
Recall that $Z_n= \left\{ q\in W |\; \gamma_n(q)=0 \right\}$ where $W$ is the neighbourhood of $L^2_0$ in $L^2_{0,\mathbb{C}}$ of Theorem \ref{KdVtheorem}. For $q\in W \setminus Z_n$ denote, as in \cite{KaP}, 
by $\theta_n(q)$ the n'th angle variable \[
\theta_n(q) = \beta_{n,n}(q) + \beta_n(q) \quad \text{and} \quad
 \beta_n(q)= \sum_{k\neq n}\beta_{n,k}(q)  \]
where \begin{align}\label{20bis}\beta_{n,n}(q)= \int_{\lambda_{2n-1}}^{\mu_n^* }\frac{\psi_n(\lambda)}{\sqrt{\Delta(\lambda)^2-4}}d\lambda \quad\text{mod} 2\pi,\end{align}
and, for $k\neq n$, 
\begin{align}\label{29ter}\beta_{n,k}(q)= \int_{\lambda_{2k-1}}^{\mu_k^* }\frac{\psi_n(\lambda)}{\sqrt{\Delta(\lambda)^2-4}}d\lambda.
\end{align}
Here $\mu_n^*$ is the point $(\mu_n, \sqrt[*]{\Delta^2(\mu_n)-4})$ on the affine curve $\Sigma_q$,
\begin{align}\label{26bis}\Sigma_q = \left\{(\lambda,z) 
\in \mathbb{C}^2 | z^2= \Delta^2(\lambda) -4\right\},\end{align}
with $\mu_n=\mu_n(q)$ denoting the n'th Dirichlet eigenvalue of the operator $-d_x^2 +q$ and $\sqrt[*]{\Delta^2(\mu_n)-4}$ denoting the square root of $\Delta^2(\mu_n)-4$ given by
\begin{align}\label{29bis}\sqrt[*]{\Delta^2(\mu_n)-4}= y_1(1,\mu_n)-y_2'(1,\mu_n).\end{align} 
The integral in \eqref{20bis}  is a straight line integral from $(\lambda_{2n-1},0)$ to $\mu_n^{*}$ in $\Sigma_q$ and the one in \eqref{29ter} is defined similarly.
See section 6 in  \cite{KaP} for more explanations concerning the notation.
Let us begin by analyzing $\beta_{n,k}$ in more detail. In Lemma 8.2 and Lemma 8.3 of \cite{KaP}, it is shown that $\beta_{n,k}$ $(k\neq n)$ is a  well defined analytic function on $W$. In  Theorem 8.5 of \cite{KaP},
 it is shown that $\beta_n= \sum_{k\neq 0} \beta_{n,k}$ is absolutely summable on $W$ and with the help of Lemma 8.4 in \cite{KaP} one proves that the following estimate holds.
\begin{lem1}\label{betaestimate} Locally uniformly on $W$, \begin{align}\label{a39bis} \beta_n= O\left(\frac{1}{n}\right).\end{align}
On $L^2_0$, \eqref{a39bis} holds uniformly on bounded subsets of $L^2_0$.
\end{lem1} 
\noindent Next let us consider $\beta_{n,n}$ in more detail. 
Recall that $\beta_{n,n}$ is defined on $W \setminus Z_n$ mod $2\pi$  and is analytic on 
$W\setminus Z_n$ mod $\pi$. Similarly as in \eqref{27bbis}, we write for $\lambda$ near $G_n$ 
\begin{align}\label{5.1} \frac{\psi_n(\lambda)}{\sqrt[c]{\Delta(\lambda)^2-4}}= \frac{\zeta_n(\lambda)}{\sqrt[s]{(\lambda_{2n}-\lambda)(\lambda -\lambda_{2n-1})}}\end{align}where \begin{align}
 \label{5.2} \zeta_n(\lambda)= (-1)^{n-1}\frac{\pi n}{\sqrt[+]{\lambda- \lambda_0}} \prod_{m\neq n} \frac{\sigma^n_m-\lambda}{\sqrt[+]{\left(\lambda_{2m}-\lambda\right)\left(\lambda_{2m-1}-\lambda\right)}}.
\end{align}   The following results improve on the asymptotic estimates of Lemma 9.2 in \cite{KaP}.
\begin{lem1}\label{xiestimate}\begin{itemize}
\item[(i)] Uniformly for  $\lambda \in U_n$ \[\zeta_n(\lambda)= 1+ \frac{1}{n}\ell^2_n.\]
\item[(ii)] Uniformly for  $\mu \in G_n$ 
\[\zeta_n(\mu)= 1+ \frac{\gamma_n}{n}\ell^2_n.\]
\item[(iii)] The error estimates in (i) and (ii) hold locally uniformly on $W$ and are uniformly bounded on bounded subsets of $L^2_0$.
\end{itemize} 

\end{lem1}
\begin{proof}
(i)
Note that by \cite{KST1}, Proposition 2.1
\begin{align}\label{27bis} 
(-1)^{n-1}\prod_{m\neq n} \frac{\sigma^n_m-\lambda}{\sqrt[+]{\left(\lambda_{2m}-\lambda\right)\left(\lambda_{2m-1}-\lambda\right)}}= 1+ \frac{1}{n}\ell^2_n 
\end{align} 
uniformly for $\lambda \in U_n$. With \[\lambda- \lambda_0= n^2\pi^2\left(1+\frac{\lambda-n^2\pi^2-\lambda_0}{n^2\pi^2}\right)\] one gets
\begin{align}\label{27ter}\frac{\sqrt[+]{\lambda-\lambda_0}}{n\pi}= \sqrt[+]{1+\frac{\lambda-n^2\pi^2-\lambda_0}{n^2\pi^2}}= 1+O\left(\frac{1}{n^2}\right).\end{align} 
Together, \eqref{27bis} and \eqref{27ter} yield 
\begin{align}\label{40bbis}
\zeta_n(\lambda)= 1 + \frac{1}{n}\ell^2_n
\end{align} 
uniformly for $\lambda\in U_n$. 

(ii) Assume that $\gamma_n\neq 0.$
 By the normalisation of the $\psi_n$-function, 
\[\int_{\lambda_{2n-1}}^{\lambda_{2n}} \frac{\psi_n(\lambda)}{\sqrt[c]{\Delta(\lambda)^2-4}}d\lambda=\pi
\] for the line integral from $\lambda_{2n-1}$ to $\lambda_{2n}$ obtained by deforming $\Gamma_n$ to the interval $[\lambda_{2n-1},\lambda_{2n}]$.
By \eqref{5.1} one then gets for any $\mu \in U_n$ 
\begin{align*}\pi&= \int_{\lambda_{2n-1}}^{\lambda_{2n}} \frac{\zeta_n(\lambda)}{\sqrt[s]{(\lambda_{2n}-\lambda)(\lambda -\lambda_{2n-1})}}d\lambda \\ &= \int_{\lambda_{2n-1}}^{\lambda_{2n}} 
\frac{\zeta_n(\mu)}{\sqrt[s]{(\lambda_{2n}-\lambda)(\lambda -\lambda_{2n-1})}}d\lambda +\int_{\lambda_{2n-1}}^{\lambda_{2n}} \frac{\zeta_n(\lambda)-\zeta_n(\mu)}{\sqrt[s]{(\lambda_{2n}-\lambda)(\lambda -\lambda_{2n-1})}}d\lambda.
\end{align*}
By a straightforward computation,
\begin{align}\label{38bbis}\int_{\lambda_{2n-1}}^{\lambda_{2n}} \frac{d\lambda}{\sqrt[s]{(\lambda_{2n}-\lambda)(\lambda -\lambda_{2n-1})}} =\pi.
\end{align} 
Combined with Lemma M.1 in \cite{KaP}, we then obtain
\begin{align}\nonumber|\zeta_n(\mu)-1| & = \frac{1}{\pi}\left|\int_{\lambda_{2n-1}}^{\lambda_{2n}} \frac{\zeta_n(\lambda)-\zeta_n(\mu)}{\sqrt[s]{(\lambda_{2n}-
\lambda)(\lambda -\lambda_{2n-1})}}d\lambda\right| \\\label{38tter}& \leq \max_{\lambda \in G_n}|\zeta_n(\lambda)-\zeta_n(\mu)|.
\end{align}
This bound holds no matter if $\gamma_n\neq 0$ or not.  Recall that $\zeta_n(\lambda)$ is analytic for $\lambda$ in $U_n$. Hence  one concludes from \eqref{40bbis} and Cauchy's estimate that
\[\dot\zeta_n(\mu)=\frac{1}{n}\ell^2_n\]  uniformly for $\mu \in U_n$ and thus by the mean value theorem, for $q\in W$ and $\mu \in G_n$,
\[\max_{\lambda \in G_n}\left|\zeta_n(\lambda)-\zeta_n(\mu)\right|\leq \frac{|\gamma_n|}{n}\ell^2_n\] uniformly for $\mu \in G_n$ as claimed. This combined with \eqref{38tter} shows  (ii). 

(iii) Going through the arguments of the proofs of (i) and (ii) one sees that the estimates hold locally uniformly on $W$ and uniformly on bounded subsets of $L^2_0$.

\end{proof}

Instead of describing the improved asymptotics of $\beta_{n,n}$, we directly study the asymptotics of  $z_n^{\pm}$, defined on $W\setminus Z_n$ by
\begin{align}\label{25bis} 
z_n^{\pm}=\gamma_n e^{\pm i \beta_{n,n}}.
\end{align}
This is the topic of the following section.

\section{Asymptotics of $z_n^{\pm}$}\label{section4bis}
The purpose of this section is to prove sharp asymptotic estimates of $z_n^{\pm}$.
First note that it was shown in \cite{KaP} that $z_n^{\pm}$, given by \eqref{25bis}, analytically extend to all of $W$ -- see formula (9.4) in \cite{KaP}. 
\begin{prop1}\label{zasymtotics} For $N\in \mathbb{Z}_{\geq 0}$, let $W_N\subseteq W\cap H^N_{0,\mathbb{C}}$ be the neighbourhood of $
H^N_0$ given in Theorem \ref{n3.2ter}. For any $q\in W_N,$ \[ z_n^{\pm} = 2 \hat q_{\mp n}
+\frac{1}{n^{N+1}}\ell^2_n\] locally uniformly on $W\cap H_{0,\mathbb{C}}^N$.
On $H^N_0$, the error is uniformly bounded on bounded subsets of $H^N_0$.
\end{prop1} 
First we need to make some preparations. Assume that $q\in W\setminus Z_n.$ Then \begin{align}\label{X1}
z_n^\pm = \gamma_n \exp\left(\pm i \int_{\lambda_{2n-1}}^{\mu_n^*}\frac{\psi_n(\lambda)}{\sqrt[]{\Delta(\lambda)^2-4}}d\lambda\right). \end{align}
If $\mu_n=\lambda_{2n-1}$, then $z_n^\pm =\gamma_n$ whereas if $\mu_n= \lambda_{2n}$, then $z_n^\pm= -\gamma_n$ by the normalization of $\psi_n$. By Lemma 9.1 in \cite{KaP}, $z_n^\pm$ 
are analytic functions on $W\setminus Z_n.$ In the case where
$\mu_n \notin \{\lambda_{2n},\lambda_{2n-1}\}$, we choose as path of integration the interval $[\lambda_{2n-1}, \mu_n]$ in $\mathbb{C}$ and obtain the following formula 
\begin{align}\label{X2}
z_n^\pm = \gamma_n \exp\left(\pm i\int_{\lambda_{2n-1}}^{\mu_n}\frac{\psi_n(\lambda)}{\sqrt[*]{\Delta(\lambda)^2-4}}d\lambda\right) \end{align}
where for $\lambda$ in $[\lambda_{2n-1}, \mu_n],\,\sqrt[*]{\Delta(\lambda)^2-4}$ is defined to be the continuous function on $[\lambda_{2n-1},\lambda_{2n}]$ with sign determined by
\[\sqrt[*]{\Delta(\lambda)^2-4}|_{\lambda =\mu_n}=\sqrt[*]{\Delta(\mu_n)^2-4}.\]
As, by assumption, $\mu_n \notin \{\lambda_{2n},\lambda_{2n-1}\}$, the root
$\sqrt[*]{\Delta(\lambda)^2-4}$ is well-defined. In analogy with formula \eqref{5.1} define
\begin{align}\label{42bbis}\sqrt[*]{(\lambda_{2n}-\lambda)(\lambda-\lambda_{2n-1})}:= \sqrt[*]{\Delta(\lambda)^2-4}\cdot \frac{\zeta_n(\lambda)}{\psi_n(\lambda)},\qquad \lambda \in [\lambda_{2n-1},\mu_n ].\end{align} 
As $\sqrt[*]{\Delta(\mu_n)^2-4}= y_1(1,\mu_n)-y_2'(1,\mu_n)$ is defined on all of $W$ and analytic there, $\sqrt[*]{(\lambda_{2n}-\mu_n)(\mu_n-\lambda_{2n-1})}$ analytically extends to $W$ as well and
\begin{align}\label{X3}
z_n^\pm = \gamma_n \exp\left(\pm i\int_{\lambda_{2n-1}}^{\mu_n}\frac{\zeta_n(\lambda)}{\sqrt[*]{(\lambda_{2n}-\lambda)(\lambda-\lambda_{2n-1})}}d\lambda\right) .\end{align} To obtain the claimed 
estimates for $z_n^\pm$, we write $z_n^\pm$  as a product 
\begin{align}\label{X4}
z_n^+= u_n^+ v_n^+\quad \text{and} \quad z_n^-= u_n^- v_n^- \end{align}
where 
\begin{align}\label{X5}
u_n^\pm = \gamma_n \exp\left(\pm i\int_{\lambda_{2n-1}}^{\mu_n}\frac{1}{\sqrt[*]{(\lambda_{2n}-\lambda)(\lambda-\lambda_{2n-1})}}d\lambda\right) \end{align}
and
\begin{align}\label{X6}
v_n^\pm =  \exp\left(\pm i\int_{\lambda_{2n-1}}^{\mu_n}\frac{\zeta_n(\lambda)-1}{\sqrt[*]{(\lambda_{2n}-\lambda)(\lambda-\lambda_{2n-1})}}d\lambda\right) \end{align}
Arguing as in the proof of Lemma 9.1 in \cite{KaP} one concludes that $u_n^\pm$ are analytic on $W\setminus Z_n$. Together with the analyticity of $z_n^\pm$ on $W\setminus Z_n$ it then follows that 
$v_n^\pm$  are analytic on $W\setminus Z_n$ as well.
For  $q\in W \setminus Z_n$ with $\mu_n =\lambda_{2n-1}$, one easily sees that
\begin{align}\label{X7}
u_n^\pm= \gamma_n \quad \text{and} \quad v_n^\pm= 1.\end{align}
A straightforward computation shows that
for $q\in W \setminus Z_n$ with $\mu_n =\lambda_{2n}$ 
\begin{align}\label{X.10}
u_n^\pm = -\gamma_n \quad v_n^\pm = 1. \end{align}
 We will see that both, $u_n^\pm$ and $v_n^\pm$, continuously extend to all of $W$ and that they admit asymptotics for $n\rightarrow \infty$ which allow to prove the asymptotics of $z_n^\pm$, claimed in Proposition \ref{zasymtotics}.
The quantities $u_n^\pm$ and $v_n^\pm$ will be studied separately. Let us begin with the $u_n^\pm$'s. To this end introduce for any $q$ in $W$,
\begin{align}\label{X.11}\Delta_n(\lambda)= 4\,\frac{\lambda -\lambda_0}{\pi_n^2}\left(\prod_{m\neq n}\frac{\lambda_{2m}-\lambda}{\pi_m^2}\right)\left(\prod_{m\neq n}\frac{\lambda_{2m-1}-\lambda}{\pi_m^2}\right).\end{align}
Hence for $\lambda$ in $U_n$, the principal branch of the square root $\sqrt[+]{\Delta_n(\lambda)}$ is well-defined. Furthermore recall that in section 4, we have introduced for any $n\geq 1$ and $q\in W$
\begin{align}\label{X.20} \kappa_n(q)= \log{(-1)^n y_2'(1,\mu_n)}
\end{align} where $\mu_n=\mu_n(q)$ is the n'th Dirichlet eigenvalue. Here log denotes the principal branch of the logarithm.
\begin{prop1}\label{lemX.10} For any $n\geq 1$ and $q\in W \setminus Z_n$ 
\begin{align}\label{X.12}u_n^\pm= 2(\tau_n- \mu_n)\pm i\frac{2\pi n}{\sqrt[+]{\Delta_n(\mu_n)}}2\sinh \kappa_n.\end{align}
In particular, for any $n\geq 1$, $u_n^{\pm}$ extend analytically to all of $W$. For $q\in Z_n$, \begin{align}\label{X.16b} u_n^\pm= 2(\tau_n -\mu_n)\pm 2i \sqrt[*]{-(\tau_n-\mu_n)^2}.
\end{align}
\end{prop1}
\begin{rem1}\label{remark4bis.1}
In Appendix A, Proposition \ref{lemX.10} is used to derive the formula for the differential of the Birkhoff map at $q=0$ by a short calculation.
\end{rem1}
\begin{proof}[Proof of Proposition \ref{lemX.10}]By the definition \eqref{X.20} of $\kappa_n$, \[2 \sinh \kappa_n = (-1)^n y_2'(1,\mu_n) -\left((-1)^n y_2'(1,\mu_n)\right)^{-1}= (-1)^{n-1}\sqrt[*]{\Delta(\mu_n)^2-4}.\] 
Recall that by \eqref{18bis},
\begin{align}\Delta^2(\lambda)-4= \frac{(\lambda_{2n}-\lambda)(\lambda- \lambda_{2n-1})}{\pi_n^2}\Delta_n(\lambda).\end{align}
By the definition \eqref{5.2} of $\zeta_n(\lambda)$ and the definition \eqref{20bis} of $\psi_n$ one sees that for $q$ real , $\zeta_n(\mu_n)>0$ and $(-1)^{n-1}\psi_n(\mu_n)>0.$ Hence by
 the definition \eqref{42bbis} of the root $\sqrt[*]{(\lambda_{2n}-\lambda)(\lambda- \lambda_{2n-1})}$ it follows that for $q$ real
\[(-1)^{n-1}\sqrt[*]{\Delta(\mu_n)^2-4}=\sqrt[+]{\Delta_n(\mu_n)}\frac{\sqrt[*]{(\lambda_{2n}-\mu_n)(\mu_n- \lambda_{2n-1})}}{\pi_n} .\]
As both sides of the last identity are analytic on $W$ it holds for any $q\in W$.
Hence it is to prove that 
\begin{align}\label{X.13}u_n^\pm= 2(\tau_n- \mu_n)\pm 2 i\sqrt[*]{(\lambda_{2n}-\mu_n)(\mu_n- \lambda_{2n-1})}.\end{align}
In the case where $\mu_n = \lambda_{2n-1}$, one obtains from \eqref{X7} that the left and right hand side of \eqref{X.13} are equal to $\gamma_n$. If $\mu_n= \lambda_{2n}$, by \eqref{X.10}, both sides of 
\eqref{X.13} equal $-\gamma_n$.  It remains to verify \eqref{X.13} for $q\in W \setminus Z_n$ with $\mu_n \notin \{\lambda_{2n},\lambda_{2n-1}\}$. Without loss of generality we may assume that
$|\mu_n-\lambda_{2n-1}|\leq |\mu_n -\lambda_{2n}|$. Otherwise, we exchange the role of $\lambda_{2n-1}$ and $\lambda_{2n}$ and note that this will not change the value of $u_n^\pm$. 
Denote by $V_n\subseteq U_n$ a (small) open neighbourhood of $(\lambda_{2n-1},\mu_n]$ which does not contain $\lambda_{2n-1}$ nor $\lambda_{2n}$.
There, the root $\sqrt[*]{(\lambda_{2n}-\mu)(\mu- \lambda_{2n-1})}$, defined on $[\lambda_{2n-1},\mu_n]$, continuously extends to $V_n\cup \{\lambda_{2n-1}\}$. We again denote it by 
$\sqrt[*]{(\lambda_{2n}-\mu)(\mu- \lambda_{2n-1})}$. Furthermore, for $\mu\in V_n$, introduce
\begin{align}\label{X.14}f_n^\pm(\mu)= \gamma_n\exp\left( \pm i\int_{\lambda_{2n-1}}^{\mu}\frac{d\lambda}{\sqrt[*]{(\lambda_{2n}-\lambda)(\lambda-\lambda_{2n-1})}}\right)\end{align} and
\begin{align}\label{X.15}g_n^\pm(\mu)= 2(\tau_n- \mu)\pm 2 i\sqrt[*]{(\lambda_{2n}-\mu)(\mu- \lambda_{2n-1})}.\end{align}
Then \[\lim_{\mu \to \lambda_{2n-1}}f_n^\pm(\mu)=\gamma_n = \lim_{\mu \to \lambda_{2n-1}}g_n^\pm(\mu), \] 
\begin{align}\nonumber \partial_\mu f_n^\pm(\mu)= \pm f_n^\pm(\mu) \frac{i}{\sqrt[*]{(\lambda_{2n}-\mu)(\mu-\lambda_{2n-1})}}\end{align}
and
\begin{align}\nonumber \partial_{\mu} g_n^\pm(\mu)= \pm g_n^\pm(\mu) \frac{i}{\sqrt[*]{(\lambda_{2n}-\mu)(\mu-\lambda_{2n-1})}}.\end{align}
Hence $f_n^+$ and $g_n^+$ satisfy the same 1st order differential equation and have the same value at $\mu=\lambda_{2n-1}$. Hence 
\begin{align}\label{X.16} f_n^+(\mu )= 
g_n^+(\mu)\quad \forall \mu \in V_n. \end{align}
In particular $f_n^+(\mu_n )= 
g_n^+(\mu_n)$. Similarly one has $f_n^-(\mu_n )= 
g_n^-(\mu_n)$ and the identity \eqref{X.13} is proved.

Note that the right hand side of \eqref{X.12} is defined on all of $W$. 
By the analyticity of $y_j(1,\lambda, q)\; (j= 1,2)$ on $\mathbb{C}\times W$, the analyticity of $\tau_n,\mu_n,\kappa_n$ on $W$, and the analyticity of $\Delta_n(\lambda,q)$ on $U_n\times W$, 
it then follows that the right hand side of \eqref{X.12} is analytic on $W$. Formula \eqref{X.16b} is obtained from \eqref{X.13}. 
\end{proof}

As an application of Proposition \ref{lemX.10} and the asymptotics of section 3 and section 4 we obtain 
\begin{cor1}\label{corX.13} For $q$ in $W_N$ with $W_N \subseteq W\cap H^N_{0,\mathbb{C}}$ given as in    Corollary \ref{corollary3.2},
\[u_n^\pm = 2\langle q,\cos{2\pi n x}\rangle\pm 2i \langle q,\sin{2\pi n x}\rangle + \frac{1}{n^{N+1}}\ell^2_n\] locally uniformly on $W_N$. On $H^N_0$, the error is uniformly bounded on bounded subsets of $H^N_0$.
\end{cor1}
\begin{proof}[Proof of Corollary \ref{corX.13}]
By Proposition \ref{lemX.10}, for $q\in W$,
\[ u_n^\pm= 2(\tau_n- \mu_n)\pm i\frac{2\pi n}{\sqrt[+]{\Delta_n(\mu_n)}}2 \sinh \kappa_n.\]  By Theorem \ref{4bis.1}, for $q$ in $W\cap H^N_{0,\mathbb{C}}$,
\[\kappa_n= \frac{1}{2\pi n}\left(\langle q,\sin{2\pi n x}\rangle+ \frac{1}{n^{N+1}}\ell^2_n \right)\]
and the Taylor expansion of $\sinh{z}$ at $z=0$ then yields
\[\sinh{\kappa_n} = \frac{1}{2\pi n}\left(\langle q,\sin{2\pi n x}\rangle+ \frac{1}{n^{N+1}}\ell^2_n \right).\]
According to \cite{KST1}, Corollary 7.1, $\Delta_n(\lambda)= 1+ \frac{1}{n}\ell^2_n$ and hence 
\[\left(\sqrt[+]{\Delta_n(\lambda)}\right)^{-1}= 1 +\frac{1}{n}\ell^2_n\]  uniformly for $\lambda \in U_n$. Furthermore, by Corollary \ref{corollary3.2}, for $q\in W_N$, 
\begin{align*}2(\tau_n-\mu_n)= 2\langle q,\cos2\pi nx\rangle +\frac{1}{n^{N+1}}\ell^2_n . \end{align*} 
Altogether we get for $q$ in $W_N$
\begin{align*}u_n^\pm= & 2\langle q,\cos2\pi nx\rangle +\frac{1}{n^{N+1}}\ell^2_n \pm 2i\left(1+ \frac{1}{n}\ell^2_n\right)\left(\langle q,\sin 2\pi nx\rangle+ \frac{1}{n^{N+1}}\ell^2_n \right)\\ 
=&2\langle q,\cos2\pi nx\rangle\pm 2i\langle q,\sin 2\pi nx\rangle + \frac{1}{n^{N+1}}\ell^2_n \end{align*} as claimed. 
Going through the arguments of the proof one sees that the estimate holds locally uniformly on $W_N$ and uniformly on bounded subsets of $H^N_0$.

\end{proof}
It remains to analyze the asymptotics for $v_n^\pm$. We need the following auxiliary result.
\begin{lem1}\label{lemX.14} For $q$ in $W\setminus Z_n$ with $|\mu_n -\lambda_{2n-1}| \leq |\mu_n - \lambda_{2n}|$
\begin{align}\label{X.21} \left|\int_{\lambda_{2n-1}}^{\mu_n}\frac{\zeta_n(\lambda)-\zeta_n(\lambda_{2n-1})}{\sqrt[*]{(\lambda_{2n}-\lambda)(\lambda-\lambda_{2n-1})}}
d\lambda \right| \leq \left(|\mu_n- \tau_n| +|\gamma_n|\right)
\sup_{\lambda \in [\lambda_{2n-1},\mu_n]}|\dot \zeta_n(\lambda)|\end{align}
and 
\begin{align}\label{X.22} \left|\int_{\lambda_{2n-1}}^{\mu_n}\frac{\zeta_n(\lambda_{2n-1})-1}{\sqrt[*]{(\lambda_{2n}-\lambda)(\lambda-\lambda_{2n-1})}}
d\lambda \right| \leq 3 \left(|\mu_n- \tau_n| +|\gamma_n|\right)
\sup_{\lambda \in G_n}\frac{
|\zeta_n(\lambda)-1|}{|\gamma_n|} \end{align} If $|\mu_n -\lambda_{2n-1}| \geq |\mu_n - \lambda_{2n}|$, \eqref{X.21} and \eqref{X.22}  hold if the roles of $\lambda_{2n-1}$ and $\lambda_{2n}$ are interchanged.
\end{lem1}
\begin{proof}[Proof of Lemma \ref{lemX.14}]
First, assume that $q\in W\setminus Z_n$ and 
\begin{align}\label{X.23}|\mu_n -\lambda_{2n-1}|\leq |\mu_n -\lambda_{2n}|.\end{align}
This implies that $\mu_n\neq \lambda_{2n}.$ If $\mu_n=\lambda_{2n-1}$ then \eqref{X.21} and \eqref{X.22} hold as their left-hand sides vanish. Now assume that $\mu_n\neq \lambda_{2n-1}$.
By the mean value theorem
\begin{align}\label{62bbis}\zeta_n(\lambda)= f_n(\lambda)(\lambda- \lambda_{2n-1}) + \zeta_n(\lambda_{2n-1})\end{align}
where for $\lambda$ in $U_n$ 
\[f_n(\lambda) = \int_0^1 \dot\zeta_n\left(\lambda_{2n-1}+s(\lambda -\lambda_{2n-1})\right)ds.\] Hence 
\begin{align}\label{62tter} \sup_{\lambda \in [\lambda_{2n-1}, \mu_n]}\left|f_n(\lambda)\right|\leq  \sup_{\lambda \in [\lambda_{2n-1}, \mu_n]}\left|\dot\zeta_n(\lambda)\right|.
\end{align}
(Here we assume (without loss of generality) that $U_n$ is convex.)  It follows from \eqref{62bbis} that, with $\lambda(t)= \lambda_{2n-1} +t(\mu_n -\lambda_{2n-1})$,
\begin{align*} \left|\int_{\lambda_{2n-1}}^{\mu_n}\frac{\zeta_n(\lambda)-\zeta_n(\lambda_{2n-1})}{\sqrt[*]{(\lambda_{2n}-\lambda)(\lambda-\lambda_{2n-1})}}
d\lambda \right| = \left| \int_0^1 f_n(\lambda) \frac{\sqrt{\lambda-\lambda_{2n-1}}}{\sqrt{\lambda_{2n}-\lambda}}(\mu_n-\lambda_{2n-1})dt\right|. \end{align*}
In view of \eqref{X.23}, $|\lambda(t)-\lambda_{2n-1}|\leq |\lambda(t)- \lambda_{2n}|$ for any $0\leq t\leq 1$ and thus
\begin{align*} &\left|\int_{0}^{1}\frac{f_n(\lambda)\sqrt{\lambda-\lambda_{2n-1}}}{\sqrt{\lambda_{2n}-\lambda}}(\mu_n-\lambda_{2n-1})
dt \right|\\ \leq & |\mu_n -\lambda_{2n-1}|\sup_{\lambda \in [\lambda_{2n-1},\mu_n]}|f_n(\lambda)| .\end{align*}
As $|\mu_n -\lambda_{2n-1}|\leq |\mu_n-\tau_n|+ |\gamma_n|/2$, the claimed estimate \eqref{X.21} then follows from \eqref{62tter}. Next consider the term \[(\zeta_n(\lambda_{2n-1})-1)\int_{\lambda_{2n-1}}^{\mu_n}\frac{d\lambda}{\sqrt[*]{(\lambda_{2n}-\lambda)(\lambda-\lambda_{2n-1})}}.
\] To estimate the integral, let
$\lambda(t)= \lambda_{2n-1} +t(\mu_n -\lambda_{2n-1})$ to get 
\begin{align}\nonumber \left|\int_{\lambda_{2n-1}}^{\mu_n}\frac{d\lambda}{\sqrt[*]{(\lambda_{2n}-\lambda)(\lambda-\lambda_{2n-1})}}\right|= & \left| \int_0^1\frac{\sqrt{\mu_n-\lambda_{2n-1}}}{\sqrt{\lambda_{2n}-\lambda}}\frac{d t}{\sqrt{t}}\right|\\\label{62ter+} \leq & \sqrt{2}\frac{2|\mu_n-\lambda_{2n-1}|^{\frac{1}{2}}}{|\gamma_n|^{\frac{1}{2}}}\end{align} where we used again that by \eqref{X.23}, $|\lambda_{2n}-\lambda(t)|\geq |\gamma_n/2|$ for $0\leq t\leq 1$. In view of Lemma \ref{xiestimate} (ii) and as 
\[2|\mu_n-\lambda_{2n-1}|^{\frac{1}{2}}|\gamma_n|^{\frac{1}{2}} \leq |\mu_n-\lambda_{2n-1}| + |\gamma_n|\leq |\mu_n-\tau_n| + \frac{3}{2}|\gamma_n|\] the claimed estimate \eqref{X.22} follows. The case when $|\mu_n -\lambda_{2n-1}| \geq |\mu_n - \lambda_{2n}|$ is treated in a similar way.
\end{proof}
\begin{cor1}\label{corX.15} For any $n\geq 1$, $v_n^+$ and $v_n^-$ continuously extend to all of $W$. These extensions are again denoted by $v_n^\pm$. For $q\in Z_n$ with $\mu_n \neq \tau_n$
\begin{align}\label{X.25}v_n^\pm = \exp \left(\pm i \int_{\tau_n}^{\mu_n}\frac{\zeta_n(\lambda)-\zeta_n(\tau_n)}{\sqrt[*]{-(\tau_n-\lambda)^2}}d\lambda\right) \end{align} 
whereas for $q \in Z_n$ with $\mu_n = \tau_n ,\,v_n^\pm =1$.
\end{cor1}
\begin{proof}[Proof of Corollary \ref{corX.15}] Let $q\in Z_n$ with $\mu_n \neq \tau_n$. It follows from \eqref{3.10}, \eqref{5.1}, and \eqref{38bbis} that for $q\in W\setminus Z_n,$
\begin{align}\label{63bbis}\int_{\lambda_{2n-1}}^{\lambda_{2n}}\frac{\zeta_n(\lambda)-1}{\sqrt[]{(\lambda_{2n}-\lambda)(\lambda-\lambda_{2n-1})}}
d\lambda=0
\end{align} for any choice of the sign of the root. In particular, \eqref{63bbis} holds for the $*$-root and thus 
\[\int_{\lambda_{2n-1}}^{\mu_n}\frac{\zeta_n(\lambda)-1}{\sqrt[*]{(\lambda_{2n}-\lambda)(\lambda-\lambda_{2n-1})}}
d\lambda=\int_{\lambda_{2n}}^{\mu_n}\frac{\zeta_n(\lambda)-1}{\sqrt[*]{(\lambda_{2n}-\lambda)(\lambda-\lambda_{2n-1})}}
d\lambda.\]
  Hence without loss of generality we may assume that $|\mu_n-\lambda_{2n-1}|\leq |\mu_n -\lambda_{2n}|$ for $q$, as otherwise we simply interchange the role of $\lambda_{2n}$ and $\lambda_{2n-1}$ in \eqref{X6}.  
Moreover, we may assume that the isolating neighbourhood $U_n$ of $q$ is also an isolating neighbourhood of $p$ for any $p$ in some neighbourhood $W_q$ of $q$. To compute the limit 
$\lim_{p\rightarrow q,\, p\in W_q\setminus Z_n}v_n^\pm$ split
\[\int_{\lambda_{2n-1}}^{\mu_n}\frac{\zeta_n(\lambda)-1}{\sqrt[*]{(\lambda_{2n}-\lambda)(\lambda-\lambda_{2n-1})}}
d\lambda\] into two parts by writing 
\[\zeta_n(\lambda)-1=(\zeta_n(\lambda)-\zeta_n(\lambda_{2n-1}))+ (\zeta_n(\lambda_{2n-1})-1).\]
Then \[\lim_{p\rightarrow q,\, p\in W_q\setminus Z_n}\int_{\lambda_{2n-1}}^{\mu_n}\frac{\zeta_n(\lambda)-\zeta_n(\lambda_{2n-1})}{\sqrt[*]{(\lambda_{2n}-\lambda)(\lambda-\lambda_{2n-1})}}
d\lambda =\int_{\tau_n}^{\mu_n}\frac{\zeta_n(\lambda)-\zeta_n(\tau_n)}{\sqrt[*]{-(\tau_n-\lambda)^2}}d\lambda\Big|_{p=q} \]
and in view of Lemma \ref{xiestimate} and \eqref{62ter+}, 
\[\lim_{p\rightarrow q,\, p\in W_q\setminus Z_n\, |\mu_n-\lambda_{2n-1}|\leq|\lambda_{2n}-\mu_n|}\int_{\lambda_{2n-1}}^{\mu_n}\frac{\zeta_n(\lambda_{2n-1})-1}{\sqrt[*]{(\lambda_{2n}-\lambda)(\lambda-\lambda_{2n-1})}}
d\lambda =0.\]
By the same reason the integral \[\int_{\lambda_{2n}}^{\mu_n}\frac{\zeta_n(\lambda_{2n})-1}{\sqrt[*]{(\lambda_{2n}-\lambda)(\lambda-\lambda_{2n-1})}}d\lambda\] converges to zero as $p\rightarrow q$ for $p\in W_q$ with 
$|\mu_n-\lambda_{2n-1}|\geq|\lambda_{2n}-\mu_n|$.
Altogether we thus have shown that 
\[\lim_{p\rightarrow q,\, p\in W_q\setminus Z_n}v_n^\pm(p) =\exp\left(\pm i\int_{\tau_n}^{\mu_n}\frac{\zeta_n(\lambda)-\zeta_n(\tau_n)}{\sqrt[*]{-(\tau_n-\lambda)^2}}d\lambda\right) \] 
as claimed. 
Now let us consider the case where $q\in Z_n$ with $\mu_n=\tau_n$. From \eqref{X.25} one concludes that 
\[\lim_{p\rightarrow q,\, p\in W_q\cap Z_n,\; \mu_n\neq \tau_n}v_n^\pm(p) =1 \] and from Lemma \ref{lemX.14} and
 Lemma \ref{xiestimate} 
one sees that \[\lim_{p\rightarrow q,\, p\in W_q\setminus Z_n,\; \mu_n\neq \tau_n}v_n^\pm(p) =1. \]
\end{proof}
It remains to study the asymptotics of $v_n^\pm$, now defined on all of $W$.
\begin{cor1}\label{corX.16} For $q$ in $W$ \begin{align} \label{a73bis} v_n^\pm= 1+\frac{1}{n}\ell^2_n\end{align}
 locally uniformly on $W$.
On $L^2_0$, \eqref{a73bis} holds uniformly on bounded subsets of $L^2_0$.
\end{cor1}
\begin{proof}[Proof of Corollary \ref{corX.16}] We want to apply Lemma \ref{lemX.14}. First note that in view of Corollary \ref{corX.15}, the estimates of Lemma \ref{lemX.14} hold on all of $W$,  not only on $W\setminus Z_n$. 
Furthermore, by shrinking the isolating neighbourhoods we get from Lemma \ref{xiestimate} and Cauchy's estimate 
\[\sup_{\lambda \in U_n}|\dot \zeta_n(\lambda)|=\frac{1}{n}\ell^2_n.\] By Lemma \ref{xiestimate}, 
\[\sup_{\lambda \in G_n} \frac{|\zeta_n(\lambda)-1|}{|\gamma_n|}=\frac{1}{n}\ell^2_n.\]  Using that $|e^x-1|\leq |x|e^{|x|}$, the claimed estimate then follows indeed from Lemma \ref{lemX.14}.
Going through the arguments of the proof one sees that \eqref{a73bis} holds locally uniformly on $W$ and uniformly on bounded subsets of $L^2_0$.
\end{proof}
\begin{proof}[Proof of Proposition \ref{zasymtotics}] By Corollary \ref{corX.15} and Lemma \ref{lemX.14}, $u_n^\pm$ and $v_n^\pm$ are defined and continuous on $W$ for any $n\geq 1$. Hence the identities \eqref{X4} extend to all of $W$,
\begin{align}\label{feb63bis}z_n^+=u_n^+v_n^+\quad \text{and}\quad z_n^-=u_n^-v_n^-.\end{align}
By Corollary \ref{corX.13} and Corollary \ref{corX.16}, it follows that for $q$ in
 $W_N$,
 \begin{align*}z_n^\pm=&\left(2\langle q,\cos 2\pi n x\rangle \pm 2 i \langle q,\sin 2\pi n x\rangle + \frac{1}{n^{N+1}}\ell^2_n\right)\left(1+\frac{1}{n}\ell^2_n\right)\\= & 2\hat q_{\mp n} + \frac{1}{n^{N+1}}\ell^2_n .
\end{align*} Going through the arguments of the proof one sees that the estimate holds locally uniformly on $W_N$ and uniformly on bounded subsets of $H^N_0$.
This proves Proposition \ref{zasymtotics}.
\end{proof}

\section{Proof of Theorem \ref{1smoothingthm}}
In this section we prove the asymptotics of the Birkhoff map claimed in Theorem \ref{1smoothingthm} and use them to derive further properties of this map.
\begin{proof}[Proof of Theorem \ref{1smoothingthm}]
For any $N\in \mathbb{Z}_{\geq 0}$, let $W_N \subseteq W\cap H^N_{0,\mathbb{C}}$ be the neighbourhood of $H^N_0$ given by Theorem \ref {n3.2ter}. For any $q\in W$, $\Phi(q)$ is given by $\Phi(q)=(z_n(q))_{n\neq 0}$, where for any $n\geq 1$, 
\begin{align}\label{n82bis} z_{-n} =x_n+i y_n = \xi_n  \frac{z_n^+}{2} e^{i\beta_n}, \quad \text{and} \quad z_n= x_n -iy_n = \xi_n \frac{z_n^-}{2}e^{-i\beta_n}.\end{align}
 By Proposition \ref{quotientprop} and Lemma \ref{betaestimate}
\begin{align}\label{a75bis}\xi_n = \frac{1}{\sqrt{\pi n}}\left(1+\frac{1}{n}\ell^2_n\right)\quad \text{and}\quad e^{i\beta_n}= 1 + O\left(\frac{1}{n}\right) \end{align}
whereas by Proposition \ref{zasymtotics}, for $q\in W_N,$
\begin{align}\label{a75ter} \frac{z_n^+}{2} =\hat q_{-n}
+\frac{1}{n^{N+1}}\ell^2_n. \end{align}
 Hence
\[z_{-n} =\frac{1}{\sqrt{n\pi}}\left(\hat q_{-n}+ \frac{1}{n^{N+1}}\ell^2_n\right).\]
A similar estimate holds for $z_n$. Note that all the asymptotic estimates referred to hold locally uniformly on $W_N$.
As by Theorem \ref{KdVtheorem}, $z_{\pm n}$ are analytic on $W$ for any $n\geq 1$ it follows from Theorem A.5 in \cite{KaP} that $\Phi-\Phi_0:W_N\to \mathfrak{h}^{N+3/2}_\mathbb{C}$ is analytic.
This proves the first part of Theorem \ref{1smoothingthm}.

Going through the above arguments and using that by Proposition \ref{quotientprop}, Lemma \ref{betaestimate}, and Proposition \ref{zasymtotics}, the estimates \eqref{a75bis} and \eqref{a75ter} hold uniformly on bounded sets of $H^N_0$, it follows that 
 $A:=(\Phi-\Phi_0):H^N_0\to \mathfrak{h}^{N+3/2}$ is bounded. 
\end{proof}
\vspace{0.4cm}
\begin{prop1}\label{Bbounded}Let $N\in \mathbb{Z}_{\geq 0}$. The restriction of $\Phi^{-1}$  to $\mathfrak{h}^{N+1/2}$ is of the form 
$\Phi^{-1}= \Phi_0^{-1} +B $, with $B=- \Phi^{-1}_0\circ A\circ \Phi^{-1}$. $B$ is a map from $ \mathfrak{h}^{N+1/2}$ to $H^{N+1}_0 $. It is bounded and real analytic. 
\end{prop1}
\begin{proof} First let us verify the formula for $B$. Clearly, it follows from \begin{align*}
Id_{\mathfrak{h}^{N+1/2}}  = (\Phi_0 +A)\circ \Phi^{-1} = Id_{\mathfrak{h}^{N+1/2}}  +\Phi_0\circ B + A\circ \Phi^{-1}
\end{align*}
that
 \[B=- \Phi^{-1}_0\circ A\circ \Phi^{-1}.\]
As $B$ is given by the composition of real analytic maps, 
\[\mathfrak{h}^{N+1/2} \overset{\Phi^{-1}}{\longrightarrow} H^N_0 \overset{A}{\longrightarrow} \mathfrak{h}^{N+3/2} \overset{-\Phi_0^{-1}}{\longrightarrow} H^{N+1}_0,   \]
it is itself real analytic.  
It remains to prove that for any $N\in \mathbb{Z}_{\geq 0}$, $B:\mathfrak{h}^{N+1/2}\to H^{N+1}$ is bounded.
First note that for any $N\in \mathbb{Z}_{\geq 0}$, the inverse of the weighted Fourier transform 
$\Phi_0^{-1}:\mathfrak{h}^{N+1/2}\to H^N_0$ and, by Theorem \ref{1smoothingthm}, the nonlinear map 
$A: H^N_0\to \mathfrak{h}^{N+3/2} $ are bounded for any $N\in \mathbb{Z}_{\geq 0}$. Furthermore,
the boundedness of $\Phi^{-1}:\mathfrak{h}^{N+1/2}\to H^N_0$ follows, in the case $N=0$, from the identity $\sum_{n\geq 1}2\pi nI_n= \frac{1}{2}\|q\|_{L^2}$, established in \cite{KaP}, Theorem E.1, and, 
in the case $N\in \mathbb{Z}_{\geq 1}$, from \cite{KorohillII}, Theorem 2.4 and Theorem 2.6. 
More precisely, in \cite{KorohillII}  it is shown that for any $q\in H^N$ with $\lambda_0(q)=0$ and $N\in\mathbb{Z}_{\geq 1}$, $\|\partial_x^N q\|_{L^2}$  can be bounded in terms of 
$\|J_n(q)\|_{\mathfrak{h}^{N+1/2}}$. Note that for any $p\in H^N_0$, $q= p-\lambda_0(p)$ satisfies $\lambda_0(q)=0$ and hence, as $N\geq 1$,
$\|\partial_x^N q\|_{L^2}=\|\partial_x^N p\|_{L^2}$. Furthermore, the quantities $(J^2_n(q))_{n\geq 1}$, introduced in \cite{KorohillII}, formula (1.2), can be shown to coincide with the action variables
$(I_n(p))_{n\geq 1}$, introduced in \cite{KaP} (cf formula (7.2)). Combining these results it follows that
$B=- \Phi^{-1}_0\circ A\circ \Phi^{-1}: \mathfrak{h}^{N+1/2} \to H^{N+1}_0$ is bounded for any $N\in \mathbb{Z}_{\geq 0}$.
\end{proof}

\begin{proof}[Proof of Corollary \ref{coro:bnf_weak}]
Clearly, as $d_0\Phi : H^N_0\to\mathfrak{h}^{N+1/2}$ is a bounded linear map, 
it is also weakly continuous. First, assume that  $N\ge 1$ and let $(q_j)_{j\ge 1}$ be a sequence in $H^N_0$ that converges weakly in $H^N_0$ to $q\in H^N_0$.  Then by Rellich's theorem, 
$q_j\to q$ strongly in $H^{N-1}_0$. By the continuity of $A:=\Phi-d_0\Phi : H^{N-1}_0\to\mathfrak{h}^{N+1/2}$ it then follows that $A(q_j)\to A(q)$ as $j\to\infty$ strongly
in $\mathfrak{h}^{N+1/2}$. Altogether we conclude that $\Phi(q_j)$ converges to $\Phi(q)$ weakly in $\mathfrak{h}^{N+1/2}$.
For $N=0$, a slightly more complicated argument is needed. 
Assume that the sequence $(q_j)_{j\ge 1}$ in $L^2_0$ converges weakly in $L^2$ to $q\in L^2_0$. Then by Rellich's theorem, $(q_j)_{j\ge 1}$ converges strongly to $q$ in $H^{-1}_0$.
As by \cite{KMT}, $\Phi$ extends to a continuous (even real analytic) map $\Phi : H^{-1}_0\to\mathfrak{h}^{-1/2}$,
we see that $(\Phi(q_j))_{j\ge 1}$ converges strongly to $\Phi(q)$ in $\mathfrak{h}^{-1/2}$. 
In particular, for any given $n\ge 1$, $\Phi_n(q_j)\to\Phi_n(q)$ as $j\to\infty$.
In addition, as $(q_j)_{j\ge 1}$ is bounded in $L^2_0$ we see from Theorem \ref{1smoothingthm} that $(A(q_j))_{j\ge 1}$, and therefore the sequence $\Phi(q_j)=d_0\Phi(q_j)+A(q_j)$, $j\ge 1$, is bounded in $\mathfrak{h}^{1/2}$.
This implies that $\Phi(q_j)$ converges weakly to $\Phi(q)$ in $\mathfrak{h}^{1/2}$.
In a same way one proves the corresponding statement for $\Phi^{-1}$.
\end{proof}

\noindent By Cauchy's estimate, Theorem \ref{1smoothingthm} yields the following asymptotics of the differential of $\Phi$.
\begin{cor1}\label{smoothcor}
For any $q\in W_N$ with $N\geq 0$, the differential of the Birkhoff map
\[d_q \Phi: H^N_{0,\mathbb{C}} \rightarrow  \mathfrak{h}^{N+1/2}_\mathbb{C} \] satisfies the asymptotic estimate \[d_q \Phi(f)= \left(\frac{1}{\sqrt{n\pi}}\hat f_n + \mathcal{R}_n(f)\right)_{n\neq 0}\] where
 \[\left(\mathcal{R}_n(f)\right)_{n\neq 0} \in \mathfrak{h}^{N+3/2}.\]
\end{cor1}
\vspace{0.4cm}
\noindent As an immediate application of Corollary \ref{smoothcor} we obtain
\begin{cor1}\label{smoothcor2}
For any $q\in W_N$ with $N\geq 0$, 
\[d_q \Phi: H^N_{0,\mathbb{C}} \rightarrow \mathfrak{h}^{N+1/2}_\mathbb{C} \] is a compact perturbation of the (weighted) Fourier transform. In particular, it satisfies the Fredholm alternative.
\end{cor1}
\noindent Corollary \ref{smoothcor2} is already proved in \cite{KaP}. The proof argues by approximation and uses quite complicated computations. It is based on the fact that the set of finite gap potentials is 
dense in $H^N_0$. In the set-up presented here, its proof is straightforward given the asymptotic estimates stated in section 2.

\section{Proof of Theorems \ref{1smoothingthm2}, \ref{applthm5}, and Corollary \ref{coro:norms}}\label{section7bis} 
In this section we prove Theorem \ref{1smoothingthm2}, Theorem \ref{applthm5}, and Corollary \ref{coro:norms}.
First we need to derive some auxilary results. Let $N\in \mathbb{Z}_{\geq 0}$ and $c\in\mathbb{R}$ be given.
For any  $q$ in $H^N_c$, write $q=p+c$ where $p\in H^N_0$ and  $\int_0^1q(x) dx =\hat q_0 =c$. Then 
\[\mathcal{H}(p+c) = \mathcal{H}(p) + 3c \int_0^1 p^2 dx + c^3.\]
Note that $\frac{1}{2}\int_0^1 p^2 dx$ is the second Hamiltonian in the KdV hierarchy. By Parseval's identity for KdV (cf \cite{KaP}, Appendix E)
\[\frac{1}{2} \int_0^1 p^2 dx = \sum_{n\geq 1} 2\pi n I_n\]
and thus  the KdV frequencies satisfy \begin{align}\label{feb65bis}\omega_n(q)=\omega_n(p) + 12 cn\pi \quad \forall n\geq 1,
\end{align}
and the KdV flow is given by 
\begin{align}\label{j68} u(t) \equiv S^t(q) = S_c^t(p) +c
\end{align}
where $S^t_c$ denotes the flow on $H^N_0$ corresponding to the Hamiltonian 
\[\mathcal{H}_c(p):= \mathcal{H}(p)+ 3c \int_0^1 p^2 dx.\]
The equations of motion corresponding to  $\mathcal{H}_c$ read, when expressed in Birkhoff coordinates $(z_n)_{n\neq 0}$,
\[  
\dot z_n =i\omega_n^c z_n  \quad \text{and}\quad  \dot z_{-n} =-i\omega_n^c z_{-n} 
\]
where
\begin{align}\label{top84bis}\omega_n^c\equiv \omega_n^c(p)=\omega_n(p) +12 cn\pi.\end{align} 
Define $\omega_{-n}^c:=-\omega_n^c\,(n\geq 1)$.  Then 
\begin{align}\label{j69} v(t)\equiv \Omega_c^t(z)= \left(e^{i\omega_n^c t}z_n\right)_{n\neq 0}
\end{align} is the flow of $\mathcal{H}_c$ expressed in (complex) Birkhoff coordinates. Here, by a slight abuse of terminology, $\omega_n^c$ is viewed as a function of $z$. Note that
 $\Phi$ conjugates the flow maps $S^t_c$ and $\Omega_c^t$, 
\begin{align}\nonumber S^t_c =\Phi^{-1}\circ \Omega_c^t \circ \Phi.
\end{align} With  $B:= \Phi^{-1}-\Phi_0^{-1}$, we get 
\begin{align}\label{j70} S^t_c =\Phi_0^{-1}\circ \Omega_c^t \circ \Phi+B\circ \Omega_c^t \circ \Phi.
\end{align}
We now analyse the map $\Phi_0^{-1}\circ \Omega_c^t \circ \Phi$ in more detail. First note that for any $z=(z_n)_{n\neq 0}\in \mathfrak{h}^{N+ 1/2}$,
\[\Phi_0^{-1}(z)(x)=\sum_{n\neq 0}\sqrt{|n|\pi} z_ne^{2\pi i n x}.\]
Hence for any $p\in H_0^N$ and $t,c\in \mathbb{R}$, 
\begin{align}\label{j71} \Phi_0^{-1}\circ \Omega_c^t\circ \Phi(p)=\sum_{n\neq 0} e^{i\omega_n^ct}\sqrt{|n|\pi}z_n(p) e^{2\pi in x}
\end{align} where $(z_n(p))_{n\neq 0}=\Phi(p)$.
For the proof of item (i) of Theorem \ref{1smoothingthm2} we need to consider the KdV flow on all of $H^N$. For any $t\in \mathbb{R}$, introduce
\[ E^t: H^N\times \mathfrak{h}^{N+1/2}\to \mathfrak{h}^{N+1/2},\; (q,z)\mapsto (e^{i\omega_n(q)t}z_n)_{n\neq 0}.\]
Denoting by $\Pi$ the projection $\Pi:H^N \to H^N_0,\,q\mapsto q-\hat q_0$ it follows that for any
$q\in H^N$,
\begin{align}\label{top90bis}
 E^t(q,\Phi \circ \Pi(q)) = (e^{i\omega_n(q)t}z_n(\Pi(q)))_{n\neq 0} = \Omega_{\hat q_0}^t\circ \Phi(\Pi(q)).
\end{align}

\begin{lem1}\label{toplem6.0}
 For any given $N\in \mathbb{Z}_{\geq 0}$ and $t\in \mathbb{R}$, the map $E^t:H^N\times \mathfrak{h}^{N+1/2}\to \mathfrak{h}^{N+1/2}$ is continuous.
\end{lem1}
\begin{proof}
 For any $q,p\in H^N$, $z,w\in \mathfrak{h}^{N+1/2}$, and for any $K>0$ one has
\begin{align*}
 \| E^t(q,z)-E^t(p,w) \|_{N+1/2}^2= \sum_{n\neq 0} |n|^{2N+1}|e^{i\omega_n(q)t} z_n -e^{i\omega_n(p)t} w_n|^2.
\end{align*}
As the KdV frequencies are real valued functions on $H^N$
\begin{align*}
 |e^{i\omega_n(q)t} z_n -e^{i\omega_n(p)t} w_n|^2 =& |z_n-e^{i(\omega_n(p)-\omega_n(q))t}w_n|^2\\ 
\leq& 2|z_n-w_n|^2 +2 |1-e^{i(\omega_n(p)-\omega_n(q))t}|^2|w_n|^2
\end{align*} and hence
\begin{align}\nonumber
  \| E^t(q,z)-E^t(p,w) \|_{N+1/2}^2\leq  2\| z-w\|^2_{N+1/2} + 8 \sum_{|n| \geq K} |n|^{2N+1}|w_n|^2 \\\label{top88bis}
+ 2\sum_{0<|n|<K} |n|^{2N+1} |e^{i(\omega_n(p)-\omega_n(q))t}-1|^2|w_n|^2.
\end{align}
By  \eqref{feb65bis},
\[\omega_n(p)-\omega_n(q)= 12n\pi(\widehat {p-q})_0 +O(1)\]
locally uniformly on $H^N\times H^N$. Moreover, it follows from \eqref{feb65bis} and \cite{KT}, Theorem 1.9 that for any $n\neq 0$, the $n$'th frequency $\omega_n :H^N \to \mathbb{R}$ is continuous. 
This combined with \eqref{top88bis} implies the statement of the lemma. 
\end{proof}

\begin{proof}[Proof of  Theorem \ref{1smoothingthm2}]
Let $N\in \mathbb{Z}_{\geq 0}$. For any $q\in H^N$ let $p=\Pi(q)= q-\hat q_0$. Then  by \eqref{j68} and \eqref{j70}, with $c=\hat q_0$,
\[S^t(q)= S^t_c(p)+c = \Phi_0^{-1}\circ \Omega_c^t\circ \Phi(p)+ B\circ \Omega^t_c\circ \Phi(p)+c.\]
Substituting $\Phi=\Phi_0+A$ into $\Phi_0^{-1}\circ \Omega_c^t\circ \Phi(p)$ then yields
\[\Phi_0^{-1}\circ \Omega^t_c \circ \Phi(p)=\sum_{n\neq 0} e^{i\omega_n^c t}\hat p_n e^{2\pi inx}
+ \Phi_0^{-1}\circ \Omega^t_c\circ A(p)\]
where by a slight abuse of terminology we denote by $\Omega_c^t\circ A(p)$ the element
$\Omega_c^t\circ A(p)= \left(e^{i\omega_n^c(p)t}a_n(p)\right)_{n\neq 0}$ and $A(p)=\left(a_n(p)\right)_{n\neq 0}$.
Taking into account that for any $n\neq 0$, $\hat q_n=\hat p_n$ and $\omega_n^c(p) = \omega_n(q)$ we conclude
\begin{align}\nonumber  R^t(q) = & S^t(q)-\sum_{n} e^{i\omega_n t}\hat q_n  e^{2\pi inx}\\ \label{j72}  = &\;  B\circ \Omega^t_c\circ \Phi(p)+ 
 \Phi_0^{-1}\circ \Omega^t_c\circ A(p).
\end{align}Statements
(ii) and (iii) of Theorem \ref{1smoothingthm2}  then follow from Theorem \ref{1smoothingthm}, Proposition \ref{Bbounded}, and the boundedness of $\Phi= \Phi_0+A$. To prove item (i) note that by \eqref{top90bis} and \eqref{j72},
\begin{align}\label{top91bis}
 R^t(q)= (B\circ E^t)(q,\Phi \circ \Pi(q)) + (\Phi_0^{-1}\circ E^t)(q,A\circ \Pi(q)).
\end{align}
It then follows from Lemma \ref{toplem6.0}, Theorem \ref{1smoothingthm}, and Proposition \ref{Bbounded} that $R^t$ is continuous.

\noindent To prove statement (iv) write for $n\in \mathbb{Z}$ arbitrary,
\begin{align*}
  \hat u_n(t)=
 e^{i\omega_n t} (\hat q_n+\tilde\rho_n(t)),
\end{align*} where  $\omega_0:=0$ and
\[\tilde\rho_n(t) :=  \widehat R^t_n(q)e^{-i\omega_n t}.\]
By the definition of $R^t$, $\hat R_0^t(q)=0 $ implying that $\tilde \rho_0(t)$
vanishes identically. Moreover $\tilde \rho_n(0)=0$ for any $n\in \mathbb{Z}$ as $R^0(q)=0$.
 Clearly $\partial_t \hat u_n(t)=i\omega_n \hat u_n(t) + e^{i\omega_n t}\partial_t \tilde \rho_n(t)$ and when substituted into the KdV  equation 
\begin{align*}- \partial_t \hat u_n(t)+
 8i\pi^3n^3 \hat u_n(t) + 6i \pi n\sum_l \hat u_l(t) \hat u_{n-l}(t) =0,
\end{align*}
one gets
\begin{align*}
 -i\omega_n  \hat u_n(t)- e^{i\omega_n t} \partial_t \tilde\rho_n(t) + 8i\pi^3n^3 \hat u_n(t) +6 i \pi n\sum_l \hat u_l(t) \hat u_{n-l}(t) =0
\end{align*} or
\begin{align}\label{osz6.6}ie^{i\omega_n t}\partial_t \tilde\rho_n(t) =&
(\omega_n  - 8\pi^3n^3 )\hat u_n(t) - 6 \pi n\sum_{l} \hat u_l(t) \hat u_{n-l}(t).
\end{align} 
Recall from \cite{KaP}, p 229 \begin{align}\label{osz6biss}\omega_n(q)= 8\pi n(\tau_n+\lambda_0/2- \sum_{m\geq 1}(\sigma_m^n-\tau_m)),\end{align} where $\lambda_0, \tau_n$, and $\sigma_m^n$ have been introduced either in 
Section \ref{intro} or  Section \ref{section3}.
For $N\geq 1, \, H^N \hookrightarrow  L^2$ is compact. As $L^2\to \mathbb{R}, q\mapsto \lambda_0(q)$ is continuous, it then follows
 that  $H^N \to \mathbb{R}: q\mapsto\lambda_0(q)$ is compact. By  Theorem \ref{n3.2ter}
\begin{align*}\tau_n=&\hat q_0+ n^2\pi^2 + \frac{1}{n}\ell^2_n\end{align*} whereas by Proposition \ref{proposition3.5}
and Theorem \ref{nthm3.2bis}, uniformly in $n$
\begin{align*}
\sigma^n_m-\tau_m=& \frac{\gamma_m^2}{m}\ell^2_m= \frac{1}{m^{2N+1}}\ell^1_m. \end{align*}
Both asymptotic estimates hold uniformly on bounded subsets of $H^N$. Hence formula \eqref{osz6biss} leads to the asymptotics 
\begin{align}\label{osz6ter}
\omega_n(q) =8\pi^3n^3 + O(n)
\end{align}
uniformly on bounded subsets of $H^N$ and statement (iv) follows.
\end{proof}
\begin{rem1}\label{toprem6.1} 
As pointed out in  Remark \ref{top0.1}, the restrictions of $R^t$ and $\partial_t R^t$ to $H^N_c$ are real analytic. To formulate this result more precisely, for any $c\in \mathbb{R}$, denote by  $H^N_{c,\mathbb{C}}$ the complexification of $H^N_c$, 
\[H^N_{c,\mathbb{C}}:= \big\{ q\in H^N_\mathbb{C}\big|\; \int_0^1 q(x)dx=c\big\} = c+H^N_{0,\mathbb{C}}.\]
Note that $H^N_{c,\mathbb{C}}$ is \emph{not} an open complex neighbourhood of $H^N_c$ in $H^N_\mathbb{C}$ as for $q$ in $H^N_{c,\mathbb{C}}$ the value of $\hat q_0$ is kept fixed. Let $E$ be a real Banach space and denote by 
$E_{\mathbb{C}}$ its complexification. A map $f:H^N_c \to E$ is said to be real analytic if $f$ extends to an analytic map $f:W\to E_{\mathbb{C}}$ where $W$ is an open neigbourhood of $H_c^N$ in $H^N_{c,\mathbb{C}}$.
In view of the fact mentioned above that $\hat R^t_0=0$ for any $t\in\mathbb{R}$, $R^t$  maps $H^N$ into $H^{N+1}_0$ and  $\partial_t R^t$  maps $H^N$ into $H^{N-1}_0$ .
Theorem \ref{1smoothingthm2} can then be amended as follows:
\begin{itemize}
 \item[(v)] for any $N\in \mathbb{Z}_{\geq 0}$ and $c\in\mathbb{R},\; R^t_{| H^N_c}:H^N_c \to H_0^{N+1}$ is real analytic;
\item[(vi)] for any $N\in\mathbb{Z}_{\geq 1}$ and $c\in\mathbb{R}$, $\partial_t R^t_{| H^N_c}:H^N_c \to H_0^{N-1}$ is real analytic.
\end{itemize}
To see that $ R^t_{| H^N_c}$ is real analytic, note that it follows from \cite{BaeKaMi}, Theorem 1 and formula \eqref{feb65bis} that for any $q\in H^N_c$, 
\begin{align}\label{top92}
 \omega_n(q) = 8\pi^3n^3 + 12 c n\pi + O(1),
\end{align} locally uniformly in a complex neighbourhood $V$ of $H^N_c$ in $H^N_{c,\mathbb{C}}$ and that for any  $n\geq 1$, $\omega_n$ is real analytic on $V$.
 One then concludes that for any $t\in \mathbb{R}$
\[E^t : V\times \mathfrak{h}^{N+1/2}_\mathbb{C} \to \mathfrak{h}^{N+1/2}_\mathbb{C},\; (q,(z_n)_{n\neq 0}) \mapsto (e^{i\omega_n(q)t}z_n)_{n\neq 0}\] is  analytic. 
Together with the analyticity  of $\Phi$ as well as $\Phi^{-1}$ and hence of
$B=\Phi^{-1}-\Phi_0^{-1}$, it then follows from \eqref{top91bis} that $R^t_{|H^N_c}$ is real analytic. The analyticity of $\partial_t R^t_{|H^N_c}$
stated in item (vi) is proved in a similar fashion.
\end{rem1}

\begin{proof}[Proof of Corollary \ref{coro:norms}]
Let $M>0$ and $s=N+\sigma$ where $N:=[s]$ so that $0\le\sigma<1$ and $N\in\mathbb{Z}_{\ge 0}$.
By Theorem \ref{1smoothingthm2} $(iii)$ there exists $C\equiv C_{N,M}>0$ so that for any $q\in H^N$ with $\|q\|_{H^N}\le M$ and for any $t\in\mathbb{R}$,
\begin{equation}\label{eq:R^t}
\|R^t(q)\|_{H^{N+1}}\le C\,.
\end{equation}
On the other hand, for any $q\in H^s$ with $\|q\|_{H^s}\le M$ and for any $t\in\mathbb R$,
\[
\big\|\sum_{n\in\mathbb Z}e^{i\omega_n^c(q) t}{\hat q}_n e^{2\pi i n x}\big\|_{H^s}=\|q\|_{H^s}\le M\,.
\]
This together with \eqref{eq:main_identity} and  \eqref{eq:R^t} implies that for any $t\in\mathbb R$,
\[
\|S^t(q)\|_{H^s}\le \big\|\sum_{n\in\mathbb Z}e^{i\omega_n^c(q) t}{\hat q}_n e^{2\pi i n x}\big\|_{H^s}+\|R^t(q)\|_{H^s}
\le M+C\,,
\]
where we used that $\|q\|_{H^N}\le\|q\|_{H^s}$ and $\|R^t(q)\|_{H^s}\le\|R^t(q)\|_{H^{N+1}}$.
\end{proof}

Next we prove Corollary \ref{coro:kdv_weak}.

\begin{proof}[Proof of Corollary \ref{coro:kdv_weak}]
Assume that the sequence $(q_j)_{j\ge 1}$ in $H^s$ weakly converges in $H^s$ to $q\in H^s$.
Let $q_j=p_j+c_j$ where $c_j:=\int_0^1q_j(x)\,dx$ and, correspondingly, $q=p+c$.
Then $(p_j)_{j\ge 1}$ converges weakly in $H^s_0$ to $p\in H^s_0$ and $c_j\to c$ as $j\to\infty$.
Further, let $z^{(j)}:=\Phi(p_j)$ and  $z:=\Phi(p)$. Then $(z^{(j)})_{j\ge 1}\subseteq\mathfrak{h}^{s+1/2}$
and by Corollary \ref{coro:bnf_weak}, $(z^{(j)})_{j\ge 1}$ converges weakly to $z$ in $\mathfrak{h}^{s+1/2}$.
In particular, there exists $C>0$ such that 
\begin{equation}\label{eq:z^j_uniformly_bounded}
\|z^{(j)}\|_{\mathfrak{h}^{s+1/2}}\le C\quad\forall j\ge 1
\end{equation}
and for any $n\ne 0$,
\begin{equation}\label{eq:component_convergence}
z^{(j)}_n\to z_n\,\,\, \mbox{as}\,\,\, j\to\infty.
\end{equation}
Recall that for any given $t\in\mathbb{R}$, $S^t(q_j)=\Phi^{-1}(\Omega^t_{c_j} (z^{(j)}))+c_j$ where
\begin{equation}\label{eq:evolution_formula}
\Omega_{c_j}^t(z^{(j)})=\big(e^{i \omega_n(p_j)t+12\pi i n c_j t}\,z^{(j)}_n\big)_{n\ne 0}\,.
\end{equation}
It follows from Theorem 1.9 in \cite{KT} that for any $n\ne 0$, the KdV frequency $\omega_n$, when viewed as  a function of the potential,
is continuous on the Sobolev space $H^{-1}_0$. As by Rellich's theorem, $p_j\to p$ strongly in $H^{-1}_0$
it the follows that for any $n\ne 0$,
\[
\lim\limits_{j\to\infty}\omega_n(p_j)=\omega_n(p)\,.
\]
This together with \eqref{eq:component_convergence} and \eqref{eq:evolution_formula} imply that for any $n\ne 0$ the $n$-th coordinate of $\Omega_{c_j}^t(z^{(j)})$ in $\mathfrak{h}^{s+1/2}$
converges to the $n$-th coordinate of $\Omega_c^t(z)$, i.e.
\[
e^{i \omega_n(p_j)t+12\pi i n c_j t}\,z^{(j)}_n\to e^{i \omega_n(p)t+12\pi i n c t}\,z_n,
\]
as $j\to\infty$.
As for any $n\ne 0$ and $j\ge 1$, the frequency $\omega_n(p_j)$ is real-valued we see from \eqref{eq:z^j_uniformly_bounded} and \eqref{eq:evolution_formula} that
\[
\|\Omega_{c_j}^t(z^{(j)})\|_{\mathfrak{h}^{s+1/2}}=\|z^{(j)}\|_{\mathfrak{h}^{s+1/2}}\le C
\]
uniformly in $j\ge 1$. This then implies that $\Omega^t_{c_j}(z^{(j)})$ converges to $\Omega^t_c(z)$
weakly in $\mathfrak{h}^{s+1/2}$. 
Finally, Corollary \ref{coro:bnf_weak} implies that $S^t(q_j)=\Phi^{-1}(\Omega^t_{c_j} (z^{(j)}))+c_j$ converges weakly in $H^s$
to $S^t(q)=\Phi^{-1}(\Omega^t_c (z))+c$.
This completes the proof of Corollary \ref{coro:kdv_weak}.
\end{proof}

We now turn to the proof of Theorem \ref{applthm5}. 
\begin{proof}[Proof of Theorem \ref{applthm5}] Let $\varepsilon,M >0$ be given and let $q\in H^s$ with $s\in\mathbb{R}_{\ge 0}$ satisfy $\|q\|_{H^s} \leq M$. KdV is globally well-posed on $H^s$ -- see e.g. \cite{CKSTT}.
Apply $(Id-P_L)$ to $S^t(q)= \sum_{n\in \mathbb{Z}} e^{i\omega_n t}\hat q_n e^{2\pi i n x}+R^t(q) $,  
\begin{align*} 
(Id-P_L)S^t(q)= \sum_{|n| > L} e^{i\omega_n t}\hat q_n e^{2\pi i n x}+(Id-P_L)R^t(q).
\end{align*}
Note that 
\begin{align}\label{feb70bis}
\big\|\sum_{|n| > L} e^{i\omega_n t}\hat q_n e^{2\pi i n x}\big\|_{H^s}= \|(Id-P_L)q\|_{H^s}.
\end{align}
By Theorem \ref{1smoothingthm} (iii),  $B_M:=\{R^t(p) |\, t\in \mathbb{R},\, p \in H^s, \|p\|_{H^s}\leq M  \}$ is bounded in $H^{N+1}$ where $N:=[s]$ and thus by the Sobolev embedding theorem relatively compact in $H^s$.
Hence there exists $L_*\in \mathbb{N}$ such that 
for any $L\geq L_*$
\[\|(Id-P_L)r\|_{H^s}<\varepsilon\quad \forall r\in B_M.\]
Thus 
\begin{align*}
\|(Id-P_L)S^t(q)\|_{H^s} \leq & \|(Id-P_L)q\|_{H^s} + \|(Id-P_L)R^t(q)\|_{H^s}\\ \leq &
\|(Id-P_L)q\|_{H^s}+ \varepsilon
\end{align*}
and 
\[\|(Id-P_L)S^t(q)\|_{H^s}\geq 
\|(Id-P_L)q\|_{H^s}- \varepsilon.\]
\end{proof}
There are other ways of approximating the KdV flow than the one considered in Theorem \ref{applthm5}. As an alternative to  the projection of the solution of KdV onto the space of trigonometric polynomials of order $L$ one could 
involve the orthogonal projection $Q_L:\mathfrak{h}^{1/2} \to\mathfrak{h}^{1/2}$ onto the $2L$ dimensional $\mathbb{R}$-vector space \[ \{(z_k)_{k\neq 0}|\;z_k=0 \;\forall |k|>L\} \]
and study \[\Phi^{-1}\circ Q_L \circ \Phi\circ S^t_c(p)= \Phi^{-1}\circ Q_L\circ \Omega^t_c \circ \Phi(p).\] 
Results similar to the ones of Theorem \ref{applthm5}  can be obtained for such a type of approximation.

\section{Appendix A}
In this appendix we use Proposition \ref{lemX.10} to derive the formula for the differential of the Birkhoff map at $q=0$ by a short calculation -- see \cite{KaP} for an alternative, but lengthier derivation. 
 Recall from \eqref{feb63bis} and \eqref{n82bis} that 
\[z_{-n} = x_n + iy_n= \xi_n e^{i\beta_n}\frac{u_n^+}{2}v_n^+.\]
By Proposition \ref{lemX.10} 
\[\frac{u_n^+}{2}= (\tau_n-\mu_n)+ i \frac{2\pi n }{\sqrt[+]{\Delta_n(\mu_n)}}\sinh \kappa_n.\]
For $q=0,\, \tau_n-\mu_n=0,\, \kappa_n=0$, $v_n^+=1$ (Corollary \ref{corX.15}), $\xi_n= \frac{1}{\sqrt{n\pi}}$ (\cite{KaP}, Theorem 7.3), $\beta_n=0$ (\cite{KaP},  Lemma 8.4), and $\mu_n, \lambda_{2n},\lambda_{2n-1}$ 
all coincide and are equal to $\pi^2n^2$. Therefore, for $\lambda$ near $n^2\pi^2$, $\lambda\neq n^2\pi^2$, 
\begin{align*} \Delta_n \vert _{q=0} =& \frac{4\lambda}{n^2\pi^2} \left(\prod_{m\geq 1}\frac{m^2\pi^2-\lambda}{m^2\pi^2}\right)^2\left(\frac{n^2\pi^2}{n^2\pi^2-\lambda}\right)^2\\
=& \frac{4\lambda}{n^2\pi^2} \left(\frac{\sin {\sqrt{\lambda}}}{\sqrt{\lambda}}\right)^2 \frac{(n^2\pi^2)^2}{(n\pi +\sqrt{\lambda})^2(n\pi -\sqrt{\lambda})^2} \\
= & \frac{4n^2\pi^2}{(n\pi +\sqrt{\lambda})^2} \left(\frac{\sin \sqrt{\lambda}}{\sqrt \lambda - n\pi}\right)^2 .
\end{align*} 
As a consequence \[\lim _{\lambda \rightarrow n^2\pi^2} \Delta_n(\lambda)\vert_{q=0}=1.\]
Hence 
\begin{align}\nonumber\partial_q(z_{-n})\vert_{q=0} =& \frac{1}{\sqrt{n\pi}} \partial_q(u_n^+/2)\vert_{q=0}\\\label{feb80} =& \frac{1}{\sqrt{n\pi}}
\left(\partial_q(\tau_n-\mu_n) +
i 2\pi n \cosh \kappa_n\cdot \partial_q 
\kappa_n \right) \vert_{q=0}.\end{align}
Note that at $q=0$, $\kappa_n=0$ and hence $\cosh \kappa_n= 1$. Furthermore
\[\partial_q \kappa_n= \frac{1}{y_2'(1,\mu_n)}\left(\partial_q y_2'(1,\mu_n)+\dot y_2'(1,\mu_n) \cdot \partial_q\mu_n \right)\]
At $q=0$, \[y_2'(1,\mu_n)= \cos n\pi = (-1)^n\]
and \[\dot y_2'(1,\mu_n)= \partial_\lambda \partial_x \frac{\sin \sqrt{\lambda}x}{\sqrt{\lambda}} \Big\vert_{\lambda= n^2\pi^2, x=1}=0.\]
Moreover (see e.g. \cite{KaP}, p 195) 
\[\partial_q y_2'(1,\mu_n)= y_2'(1,n^2\pi^2) \cos n\pi x \frac{\sin n\pi x}{n\pi} = \frac{(-1)^n}{2n\pi}\sin 2n\pi x\] and thus
\begin{align}\label{feb81}2\pi n \cosh \kappa_n \cdot \partial_q \kappa_n \vert_{q=0} = \sin 2n\pi x.\end{align}
Next let us compute $\partial_q(\tau_n-\mu_n)$ at $q=0$.
Using Riesz projectors one has, for any $q$ near 0, 
\[\mu_n = \operatorname{Tr} \left(\frac{1}{2\pi i } \int_{\Gamma_n} \lambda (\lambda - L_{D}(q))^{-1}d\lambda\right)\] and 
\[2\tau_n = \operatorname{Tr} \left(\frac{1}{2\pi i } \int_{\Gamma_n} \lambda (\lambda - L_{p}(q))^{-1}d\lambda\right)\]
where $\Gamma_n$ is a counterclockwise oriented contour around $n^2\pi^2$ and $L_{D}(q)\; [L_{p}(q)]$ is the operator $L(q)= -d_x^2 +q$ considered on the space $H^2_{0}[0,1] \;[H^2(\mathbb{R}/2\mathbb{Z},\mathbb{R})]$. Then 
\[\langle \partial_q \mu_n, h\rangle =\operatorname{Tr} \left(\frac{1}{2\pi i } \int_{\Gamma_n} \lambda (\lambda - L_{D}(q))^{-1}h(\lambda - L_{D}(q))^{-1}d\lambda\right).\]
At $q=0$, $\sqrt{2} \sin{n\pi x}$ is the $L^2$--normalized eigenfunction corresponding  to $\mu_n=n^2\pi^2$. Hence
\begin{align*}\langle \partial_q \mu_n, h\rangle =& \langle \sqrt{2} \sin n\pi x, \frac{1}{2\pi i} \int_{\Gamma_n} \frac{\lambda}{(\lambda- n^2\pi^2)^2} h\sqrt{2}\sin n\pi x d\lambda\rangle
\\= & \int_0^1 h(x) \sin^2 n\pi x dx.
\end{align*}
Similarly at $q=0$, $\sqrt{2} \cos n\pi x$ and $\sqrt{2}\sin n\pi x$ are an orthonormal basis of the eigenspace of $L_{per}(q)$ corresponding to the periodic eigenvalue 
$\lambda_{2n}= \lambda_{2n-1}=n^2\pi^2$. Hence
\begin{align*} \langle 2\tau_n, h\rangle =& \langle \sqrt{2} \sin n\pi x, \frac{1}{2\pi i} \int_{\Gamma_n} \frac{\lambda}{(\lambda- n^2\pi^2)^2} h\sqrt{2}\sin n\pi x d\lambda\rangle \\ +&\langle \sqrt{2} \cos n\pi x, \frac{1}{2\pi i} 
\int_{\Gamma_n} \frac{\lambda}{(\lambda- n^2\pi^2)^2} h\sqrt{2}\sin n\pi x d\lambda\rangle \\ = & \int_0^1 h(x) 2(\sin^2 n\pi x + \cos^2 n\pi x) dx .
\end{align*}
As $2\sin^2 n\pi x= 1-\cos 2n\pi x$ it then follows that 
\begin{align}\label{feb82} \partial_q(\tau_n-\mu_n)= \cos 2n\pi x. \end{align}
Combining \eqref{feb80}-\eqref{feb82} then yields 
\[\partial_q z_{-n} \vert_{q=0}= \frac{1}{\sqrt{n\pi}}e^{i2\pi n x}.\]
Similar computations lead to 
\[\partial_q z_{n} \vert_{q=0}= \frac{1}{\sqrt{n\pi}}e^{-i2\pi n x}.\]

\section{Appendix B}
In \cite{ETZ}, Erdogan and Tzirakis proved the following result on the approximation of solutions of KdV on the torus by corresponding solutions of the Airy equation $v_t=L_cv$ where $L_c=-\partial_x^3 + 6 c \partial_x$.

\begin{thm1}[\cite{ETZ}]\label{ETZ}
Fix $s> -\frac{1}{2}$ and $s_1 < \min(s+1,3s+1)$. 
Consider the real valued solutions of KdV on $\mathbb{R}\times\mathbb{T}$ with initial data $u(0)=q\in H^s$. Assume that there exist $C=C(\|q\|_{H^s})>0 $ and $\alpha(s)\geq 0$ so that for any $t\in \mathbb{R}$, 
$\|u(t)\|_{H^s}\leq C(1+|t|)^{\alpha(s)}$. Then 
\[u(t)-e^{tL_c}q\in C^0(\mathbb{R},H^{s_1}) \;\text{with} \; c=\int_0^1 q(x) dx\]
and there exists $C'= C'(s,s_1, \|q\|_{H^s})> 0$ so that 
\[\|u(t)-e^{tL_c}q\|_{H^{s_1}}\leq C' (1+|t|)^{1+\alpha(s)}\quad \forall t\in \mathbb{R}.
\]
\end{thm1}
\begin{rem1}
By Corollary \ref{coro:norms}, for any $s\in\mathbb{R}_{\ge 0}$, the assumption of Theorem \ref{ETZ} on $\|u(t)\|_{H^s}$ always holds with $\alpha(s)=0$.
\end{rem1}
It turns out that in the case where $s=N\in \mathbb{Z}_{\geq 1}$, the above theorem can be derived from Theorem 1.1. In fact we prove a stronger version with $s_1=s+1$ and $\alpha(s)=0$.
\begin{thm1}\label{appthmB.2}
 For any initial data $q\in H^N$ with $N\in \mathbb{Z}_{\geq 1}$
\[\|u(t)-e^{tL_c}q\|_{H^{N+1}}\leq C (1+|t|)\quad \forall t\in \mathbb{R}
\]
where $c=[q]$ and the constant $C>0$ can be chosen uniformly for bounded subsets of $q$ in $H^N$.
\end{thm1}

\begin{rem1}
 Note that for $|t|$ large, the approximation of $u(t)$ by the corresponding solution $e^{tL_c}q$ of the Airy equation is not satisfactory. If the linear approximation is replaced by the nonlinear one of Theorem 1.1, 
involving the KdV frequencies, the $H^{N+1}$-norm of the difference of the solution and its approximation remains bounded for all time.
\end{rem1}
 
To prove Theorem \ref{appthmB.2} we need first to establish asymptotics for  the KdV frequencies. Recall that the KdV frequencies of a potential $p\in L^2_0$ are given by the formulas \eqref{osz6biss},
\[\omega_n(q)= 8\pi n\left(\tau_n+\lambda_0/2- \sum_{m\geq 1}(\sigma_m^n-\tau_m)\right)\]  where $(\sigma_m^n)_{m\neq n}$ are the zeros of the entire function
$\psi_n(\lambda)= \frac{2}{n\pi}\prod_{k\neq n} \frac{\sigma_k^n-\lambda}{\pi_k^2}$.
 By Proposition 2.2, $\sigma_m^n-\tau_m= \frac{1}{m}\gamma_m^2\ell^2_m$ and by Theorem 2.4 (i), $\tau_n= n^2\pi^2+ O(\frac{1}{n^2})$ where both estimates hold uniformly on bounded subsets of $H^1_0$. Hence we have
\begin{align}\label{appB.4} 
 \omega_n= (2\pi n)^3 + 8\pi n\left( \frac{\lambda_0}{2} - \sum_{m\leq n/2}(\sigma_m^n-\tau_m)\right)+O\left(\frac{1}{n}\right)
\end{align}
uniformly on bounded subsets of $H^1_0$. We claim that 
\begin{align}\label{appB.5} 
  \frac{\lambda_0}{2} - \sum_{m\leq n/2}(\sigma_m^n-\tau_m) = O\left(\frac{1}{n^2}\right)
\end{align} uniformly on bounded subsets of $H^1_0$. To this aim recall that the normalization factor $\frac{2}{n\pi}$ in $\psi_n(\lambda)$ is chosen is such a way that 
\[1= \frac{1}{2\pi} \int_{\Gamma_n}\frac{\psi_n(\lambda)}{\sqrt[c]{\Delta^2(\lambda)-4}}d\lambda.\]
As \[\Delta^2(\lambda)-4= -4(\lambda-\lambda_0)\prod_{k\geq 1} \frac{(\lambda_{2k}-\lambda)(\lambda_{2k-1}-\lambda)}{\pi_k^4}\]
one has in the case $\lambda_{2n}= \lambda_{2n-1}$ by Cauchy's formula
\begin{align}\label{appB.6}
1=\frac{1}{2\pi i} \int_{\Gamma_n}\frac{f_n(\lambda-\tau_n)}{\lambda-\tau_n} d\lambda= f_n(0) 
\end{align}
where, with $\mu:=\lambda-\tau_n$,
\[f_n(\mu) = \frac{n\pi}{\sqrt{\lambda-\lambda_0}} (-1)^{n-1} \prod_{k\neq n} \frac{\sigma_k^n-\lambda}{\sqrt{(\lambda_{2k}-\lambda)(\lambda_{2k-1}-\lambda)}}\]
where here and in the sequel, $\sqrt{\cdot}= \sqrt[+]{\cdot}$. In the case $\lambda_{2n-1}<\lambda_{2n}$, it turns out that the identity \eqref{appB.6} holds up to an error term.
\begin{lem1}\label{applemB.2}
 $f_n(0) = 1+ O(\frac{1}{n^4})$ uniformly on bounded subsets of potentials in $H^1_0$.
\end{lem1}
\begin{proof}[Proof of Lemma \ref{applemB.2}]
 Note that for $n\geq n_0$, $f_n$ is analytic on the isolating neighbourhood  $U_n= \{|\lambda-n^2\pi^2| < r\pi^2 \}$ where $n_0$ and $r>0$ can be chosen uniformly on bounded subsets of potentials in $L^2_0$. 
Note that if $\lambda_{2n-1}<\lambda_{2n}$, then 
\begin{align*}1=& \frac{1}{\pi} \int_{-\gamma_n/2}^{\gamma_n/2}f_n(\mu)\frac{d\mu}{\sqrt{\gamma_n^2/4-\mu^2}}   \\
=& \frac{1}{\pi} \int_{0}^{\gamma_n/2}(f_n(\mu)+f_n(-\mu))\frac{d\mu}{\sqrt{\gamma_n^2/4-\mu^2}}.
\end{align*}
Now  expand $f_n$ near $\mu=0$ to get 
\[f_n(\mu)= f_n(0)+ f_n'(0)\mu +\frac{1}{2} f_n''(\theta_+) \mu^2\] and
\[f_n(-\mu)= f_n(0)- f_n'(0)\mu +\frac{1}{2} f_n''(\theta_-) \mu^2\]
where $0<\theta_+\equiv \theta_+(\mu) < \gamma_n/2$ and $-\gamma_n/2<\theta_-\equiv \theta_-(\mu)<0$. Hence 
\[f_n(\mu)+f_n(-\mu)= 2f_n(0)+ \frac{1}{2}(f''_n(\theta_+)+ f''_n(\theta_-)) \mu^2\]
and therefore
\[1= f_n(0)\frac{2}{\pi}\int_0^{\gamma_n/2}\frac{d\mu}{\sqrt{\gamma_n^2/4-\mu^2}}+ \frac{1}{2\pi}\int_0^{\gamma_n/2}(f''_n(\theta_+)+ f''_n(\theta_-))  \frac{\mu^2d\mu}{\sqrt{\gamma_n^2/4-\mu^2}}.
\]
As $\int_0^{\gamma_n/2}\frac{d\mu}{\sqrt{\gamma_n^2/4-\mu^2}}=\frac{\pi}{2}$ it then follows that 
\begin{align}\label{appB.7}
 |f_n(0)-1|\leq \max_{|\mu|\leq \gamma_n/2}|f_n''(\mu)|\cdot (\frac{\gamma_n}{2})^2.
\end{align}
Using that $f_n(\lambda-\tau_n)$ is analytic on $U_n$ for any $n\geq n_0$ one can bound $f_n''(\mu)=\partial_\mu^2 (f_n(\mu)-1)$ by Cauchy's estimate
\begin{align}\label{appB.8}
\max_{|\mu|\leq \gamma_n/2}|f_n''(\mu)|\leq C\sup_{\lambda\in U_n}|f_n(\lambda-\tau_n)-1|.
\end{align}
Let us now estimate $\sup_{\lambda\in U_n}|f_n(\lambda-\tau_n)-1|$. 
For $\lambda\in U_n$ with $n\geq n_0$, 
\[\sqrt{\lambda-\lambda_0}= n\pi \sqrt{1-\frac{\lambda_0+n^2\pi^2-\lambda}{n^2\pi^2}}= n\pi\left(1+O(\frac{1}{n^2})\right)\]
and thus 
\begin{align}\label{appB.9}
\frac{n\pi}{\sqrt{\lambda-\lambda_0}}= 1+O\left(\frac{1}{n^2}\right)
\end{align} uniformly for $\lambda \in U_n$. 
Let $\alpha_k^n= \sigma_k^n-\tau_k$, and write for $n\geq n_0$
\begin{align}\nonumber
 (-1)^{n-1}\prod_{k\neq n}\frac{\sigma_k^n-\lambda}{\sqrt{(\tau_k-\lambda)^2-\gamma_k^2/4}}= \prod_{k\neq n}\frac{(1+\frac{\alpha_k^n}{\tau_k-\lambda})}{\sqrt{1-\left(\frac{\gamma_k/2}{\tau_k-\lambda}\right)^2}} \\ \label{appB.10} = 
\exp \left(\sum_{k\neq n}\left(\log\left(1+\frac{\alpha_k^n}{\tau_k-\lambda}\right) -\frac{1}{2} \log \left(1- \left(\frac{\gamma_k/2}{\tau_k-\lambda}\right)^2\right)\right)\right).
\end{align}
Here we have chosen $n_0$ larger, if necessary, to ensure that $|\frac{\gamma_k}{\tau_k-\lambda}|,|\frac{\alpha_k^n}{\tau_k-\lambda}|\leq \frac{1}{2}$ for any $k\neq n$ and $n\geq n_0$.
Furthermore note that
\[\sup_{\lambda\in U_n}\left|\frac{\gamma_k}{\tau_k-\lambda}\right|^2= O\left(\frac{\gamma_k^2}{(k^2-n^2)^2}\right),\quad
\sup_{\lambda\in U_n}\left|\frac{\alpha_k^n}{\tau_k-\lambda}\right|^2= O\left(\left|\frac{\alpha_k^n}{k^2-n^2}\right|\right).\]
As $|\tau_k-\lambda|\geq C|k^2-n^2|\geq C|k-n||k+n|$ it then follows that uniformly for $\lambda\in U_n$
\begin{align*}
 \sum_{k\geq 1} \left|\log \left(1- \left(\frac{\gamma_k/2}{\tau_k-\lambda}\right)^2\right)\right| \leq & C \frac{1}{n^2} \sum_{k\geq 1}\gamma_k^2
\\ \sum_{k\leq n/2} \left|\log \left(1+ \frac{\alpha_k^n}{\tau_k-\lambda}\right)\right| \leq & C \frac{1}{n^2} \sum_{k\geq 1}\gamma_k^2
\end{align*}
and 
\[ \sum_{k> n/2} \left|\log \left(1+ \frac{\alpha_k^n}{\tau_k-\lambda}\right)\right| \leq C \frac{1}{n} \sum_{k> n/2}\frac{\gamma_k^2}{k}\leq C \frac{2}{n^2} \sum_{k\geq 1}\gamma_k^2. \]
Substituting these estimates into \eqref{appB.10} one gets 
\[(-1)^{n-1}\prod_{k\neq n}\frac{\sigma_k^n-\lambda}{\sqrt{(\tau_k-\lambda)^2-\gamma_k^2/4}}=1+ O\left(\frac{1}{n^2}\right)\] uniformly for $\lambda\in U_n$ which together with \eqref{appB.9} leads to 
$\sup_{\lambda\in U_n}|f_n(\lambda-\tau_n)-1|\leq C\frac{1}{n^2}$. Hence \eqref{appB.7} yields the estimate 
\[|f_n(0)-1|\leq C\frac{1}{n^2}\cdot \frac{1}{n^2}|n\gamma_n|^2.\]
Going through the arguments of the proof one sees that $f_n(0)=1+ O(\frac{1}{n^4})$ uniformly on bounded subsets of $H^1_0$.
\end{proof}

We are now in a position to prove the claimed estimate \eqref{appB.5}.
\begin{lem1}\label{applemB.3}
 Uniformly on bounded subsets of $H^1_0$, 
\[\frac{\lambda_0}{2} - \sum_{k\leq n/2}(\sigma_k^n-\tau_k)) = O\left(\frac{1}{n^2}\right).\]
\end{lem1}
\begin{proof}[Proof of Lemma \ref{applemB.3}] The starting point is the product representation of $f_n(0)$,
 \[f_n(0)= \frac{n\pi}{\sqrt{\tau_n-\lambda_0}}(-1)^{n-1}\prod_{k\neq n}\frac{\sigma_k^n-\tau_n}{\sqrt{(\tau_k-\tau_n)^2-\gamma_k^2/4}}.\]
In contrast to the proof of Lemma \ref{applemB.2} we need to expand $f_n(0)$ to higher order. To begin we consider 
\[\sqrt{\tau_n-\lambda_0}= \sqrt{\tau_n}\left(1-\frac{\lambda_0}{\tau_n}\right)^{1/2}= \sqrt{\tau_n}\left(1-\frac{1}{2}\frac{\lambda_0}{\tau_n}+ O\left(\frac{1}{n^4}\right)\right).\] 
As $\sqrt{\tau_n}= n\pi + O(\frac{1}{n^3})$ it then follows that 
\[\sqrt{\tau_n-\lambda_0}= n\pi -\frac{1}{2}\frac{\lambda_0}{n\pi} + O\left(\frac{1}{n^3}\right)\] and hence 
\begin{align}\label{appB.12}
 \frac{n\pi}{\sqrt{\tau_n-\lambda_0}}= 
 \frac{1}{1-\frac{\lambda_0}{2n^2\pi^2}+O(\frac{1}{n^4})}= 1+\frac{\lambda_0}{2n^2\pi^2}+ O\left(\frac{1}{n^4}\right).
\end{align}
Next note that with $\alpha_k^n= \sigma_k^n-\tau_k$,
\begin{align}\nonumber
 (-1)^{n-1}\prod_{k\neq n}\frac{\sigma_k^n-\tau_n}{\sqrt{(\tau_k-\tau_n)^2-\gamma_k^2/4}}= \prod_{k\neq n}\frac{(1+\frac{\alpha_k^n}{\tau_k-\tau_n})}{\sqrt{1-\left(\frac{\gamma_k/2}{\tau_k-\tau_n}\right)^2}} \\ \nonumber = 
\exp \left(\sum_{k\neq n}\left(\log\left(1+\frac{\alpha_k^n}{\tau_k-\tau_n}\right) -\frac{1}{2} \log \left(1- \left(\frac{\gamma_k/2}{\tau_k-\tau_n}\right)^2\right)\right)\right)
\end{align}
where for the latter identity we again assumed that $n\geq n_0$ and  $n_0$ is chosen so large that $|\frac{\gamma_k}{\tau_k-\tau_n}|,|\frac{\alpha_k^n}{\tau_k-\tau_n}|\leq \frac{1}{2}$ for any $k\neq n$ and $n\geq n_0$. 
Furthermore note that $|\tau_k-\tau_n|\geq C'|k-n|(k+n)\geq Cn$ and thus
\begin{align*}
 \sum_{k\leq n/2}\left(\frac{\gamma_k/2}{\tau_k-\tau_n}\right)^2=& O\left(\frac{1}{n^4}\right)\\
\sum_{k\geq n/2}\left(\frac{\gamma_k/2}{\tau_k-\tau_n}\right)^2=& O\left(\frac{1}{n^4}\sum_{k\geq 1}(k\gamma_k)^2\right)\\
\sum_{k\leq n/2}\frac{\alpha_k^n}{\tau_k-\tau_n}= &-\frac{1}{n^2\pi^2}\sum_{k\leq n/2}\alpha_k^n + \sum_{k\leq n/2}\alpha_k^n\left(\frac{1}{\tau_k-\tau_n}+ \frac{1}{n^2\pi^2}\right) \\
= & -\frac{1}{n^2\pi^2}\sum_{k\leq n/2}\alpha_k^n + O\left(\frac{1}{n^4}\sum_{k\geq 1}(k\gamma_k)^2 \right),
\end{align*}
where for the latter identity we used that 
\[\alpha_k^n\cdot \left(\frac{1}{\tau_k-\tau_n}+ \frac{1}{n^2\pi^2}\right)= \alpha_k^n\cdot \frac{\tau_k+n^2\pi^2-\tau_n}{(\tau_k-\tau_n)n^2\pi^2}= O\left(\frac{1}{k}(k\gamma_k)^2\frac{1}{n^4}\right).\]
We thus obtain the asymptotics 
\begin{align}\label{appB.14}
 (-1)^{n-1}\prod_{k\neq n}\frac{\sigma_k^n-\tau_n}{\sqrt{(\tau_k-\tau_n)^2-\gamma_k^2/4}}= 1-\frac{1}{n^2\pi^2}\sum_{k\leq n/2} \alpha_k^n+ O\left(\frac{1}{n^4}\right).
\end{align} Combining \eqref{appB.12} and \eqref{appB.14} one then gets 
\begin{align*}
 f_n(0)= &\left(1+ \frac{1}{2}\frac{\lambda_0}{n^2\pi^2} + O(\frac{1}{n^4})\right)\left(1-\frac{1}{n^2\pi^2 }\sum_{k\leq n/2} \alpha_k^n+ O(\frac{1}{n^4})\right) \\ 
=& 1+ \frac{1}{n^2\pi^2} \left(\frac{\lambda_0}{2}-\sum_{k\leq n/2} \alpha_k^n \right) + O\left(\frac{1}{n^4}\right).
\end{align*}
As by Lemma \ref{applemB.2}, $f_n(0)= 1+ O(\frac{1}{n^4})$ the claimed estimate $\frac{\lambda_0}{2}-\sum_{k\leq n/2} \alpha_k^n= O(\frac{1}{n^2})$ follows. 
Going through the arguments of the proof one sees that the latter asymptotics hold uniformly on bounded subsets of $H^1_0$.
\end{proof}

\begin{prop1}\label{apppropB.4}
 For any $q\in H^1_0$ ,
\[\omega_n- (2n\pi)^3= O\left(\frac{1}{n}\right)\] uniformly on bounded subsets of $q$ in $H^1_0$.
\end{prop1}
\begin{proof}[Proof of Proposition \ref{apppropB.4}] In view of formula \eqref{appB.4}, the claimed statement follows from Lemma \ref{applemB.3}. 
\end{proof}

\begin{proof}[Proof of Theorem \ref{appthmB.2}] Denote by $s_n^c$ the frequencies of the Airy equation 
$v_t=L_c v$ where $L_c v =-\partial_x^3 v + 6 c \partial_x v$ on the torus $\mathbb{T}$, 
\[s_n^c= (2\pi n)^3 + 12 cn\pi \quad \forall n\in \mathbb{Z}.\]
The solution of $v_t=L_cv$ with $v(0)= q$ then is given by 
\[v(t)= \sum_{k\in \mathbb{Z}}e^{is_n^c t}\hat q_ne^{2\pi inx}.\]
 On the other hand, using the notation of Theorem \ref{1smoothingthm2} and its proof, the solution $u(t)= S^t(q)$ of the KdV equation with initial data $u(0)= q$ is given by 
\[u(t)= \sum_{k\in \mathbb{Z}} e^{i\omega_n^c t}\hat q_ne^{2\pi inx}+ R^t(q)\] where $\omega_0^c=0$ and for any $n\geq 1$ 
\[\omega_n^c= \omega_n+12cn \pi\quad \text{and} \quad \omega_{-n}^c= -\omega_n^c\]
with $c= [q]$ and $\omega_n= \omega_n(q-c),\, \omega_{-n}= -\omega_n$.
As $\omega_0^c= s_0^c=0$ it then follows that 
\[u(t)-v(t)= \sum_{n\neq 0}\left(e^{i\omega_n^c t}-e^{is_n^c t}\right)\hat q_n e^{i2\pi nx}+ R^t(q)\]
and thus for any $t\in \mathbb{R}$,
\begin{align*}
 \|u(t)-v(t)\|_{H^{N+1}}\leq \left(\sum_{n\neq 0}n^{2N+2}|e^{i\omega_n^c t}-e^{is_n^c t}|^2|\hat q_n|^2\right)^{1/2} + \|R^t(q)\|_{H^{N+1}}.
\end{align*}
By Theorem \ref{1smoothingthm2}, $R^t(q)$ has the property that $\sup_{t\in \mathbb{R}}\|R^t(q)\|_{H^{N+1}}$ is uniformly bounded on bounded subsets of potentials in $H^N$. Furthermore, using that 
\begin{align*}
 |e^{i\omega_n^c t}-e^{is_n^c t}|= & |e^{i(\omega_n^c -s_n^c) t}-1|\\ 
=& \left|i(\omega_n^c-s_n^c)\int_0^t e^{i(\omega_n^c -s_n^c) \tau}d\tau\right|\\ 
\leq & |\omega_n^c-s_n^c||t|
\end{align*}
and $\omega_n^c-s_n^c= \omega_n-(2\pi n)^3$ one concludes that 
\[ \left(\sum_{n\neq 0}n^{2N+2}|e^{i\omega_n^c t}-e^{is_n^c t}|^2|\hat q_n|^2\right)^{1/2}\leq |t|  \left(\sum_{n\neq 0}n^{2N+2}|\omega_n-(2\pi n)^3|^2|\hat q_n|^2\right)^{1/2}\]
Hence the claimed estimate follows from Proposition \ref{apppropB.4}.
\end{proof}

\end{document}